\documentclass[10pt,epsfig]{amsart}
\usepackage{amsbsy,amssymb,amscd,amsfonts,latexsym,amstext,delarray,
amsmath,epsfig}  \headsep=15pt
\input xypic

 \def\n{{\noindent}}

\newtheorem{thm}{Theorem}[section]
\newtheorem{prop}[thm]{Proposition}

\newtheorem{lem}[thm]{Lemma}
\newtheorem{defn}[thm]{Definition}

\newtheorem{ex}[thm]{Example}

\numberwithin{equation}{section}

\newcommand{\ie}{{\it i.e.\/}\ }
\newcommand{\eg}{{\it e.g.\/}\ }
\newcommand{\cf}{{\it cf.\/}\ }

\def\qqq{\,,\quad \forall}

\def\A{{\mathbb A}}

\def\C{{\mathbb C}}

\def\H{{\mathbb H}}
\def\Hb{{\mathbb H}}

\def\J{{\mathbb I}}

\def\N{{\mathbb N}}

\def\P{{\mathbb P}}
\def\Qb{{\mathbb Q}}
\def\Q{{\mathbb Q}}

\def\R{{\mathbb R}}

\def\Z{{\mathbb Z}}

\def\sA{{\mathcal A}}
\def\cA{{\mathcal A}}

\def\cE{{\mathcal E}}
\def\sE{{\mathcal E}}
\def\E{{\mathcal E}}

\def\cG{{\mathcal G}}

\def\cH{{\mathcal H}}

\def\I{{\mathcal I}}
\def\Jc{{\mathcal J}}
\def\sJ{{\mathcal J}}

\def\cK{{\mathcal K}}

\def\Oc{{\mathcal O}}
\def\O{{\mathcal O}}
\def\cO{{\mathcal O}}

\def\cR{{\mathcal R}}

\def\cT{{\mathcal T}}
\def\Uc{{\mathcal U}}
\def\cU{{\mathcal U}}

\def\n{{\mathfrak n}}

\def\a{\alpha}

\def\g{\gamma}

\def\L{\Lambda}

\def\part{\partial}

\def\qqq{\,,\quad \forall}

\def\text{\hbox}

\def\Aut{{\rm{Aut}}}

\def\End{\mathop{\rm End}\nolimits}
\def\Gal{{\rm Gal}}
\def\GL{{\rm GL}}
\def\Hom{\mathop{\rm Hom}\nolimits}

\def\M{{\rm M}}

\def\SL{{\rm SL}}

\def\Tr{{\rm Tr}}

\begin{document}

\title[KMS and CM]{KMS states and complex
multiplication}

\author[Connes]{Alain Connes}
\author[Marcolli]{Matilde Marcolli}
\author[Ramachandran]{Niranjan Ramachandran}
\address{A.~Connes: Coll\`ege de France \\
3, rue d'Ulm \\ Paris, F-75005 France} \email{alain\@@connes.org}
\address{M.~Marcolli: Max--Planck Institut f\"ur Mathematik  \\
Vivatsgasse 7 \\
Bonn, D-53111 Germany} \email{marcolli\@@mpim-bonn.mpg.de}
\address{N.~Ramachandran: Department of Mathematics \\ University of
Maryland \\ College Park \\ MD 20912 USA}
\email{atma@math.umd.edu}

\maketitle

\section{Introduction}

Several results point to deep relations between noncommutative
geometry and class field theory (\cite{BC}, \cite{CM}, \cite{LvF},
\cite{Man1}). In \cite{BC} a quantum statistical mechanical system
is exhibited, with partition function the Riemann zeta function,
and whose arithmetic properties are related to the Galois theory
of the maximal abelian extension of $\Q$. In \cite{CM}, this
system is reinterpreted in terms of the geometry of commensurable
1-dimensional $\Q$-lattices, and a generalization is constructed
for 2-dimensional $\Q$-lattices. The arithmetic properties of this
$\GL_2$-system and its KMS states at zero temperature related to
the Galois theory of the modular field. The ground states and the
Galois properties are analyzed in \cite{CM} for the generic case
of elliptic curves with transcendental $j$-invariant. As the
results of \cite{CM} show, one of the main new features of the
$\GL_2$-system is the presence of symmetries by {\em
endomorphism}, through which the full Galois group of the modular
field appears as symmetries acting on the KMS equilibrium states
of the system.

\smallskip

In both the original BC system and in the $\GL_2$-system, the
arithmetic properties of zero temperature KMS states rely on an
underlying result of compatibility between ad\`elic groups of
symmetries and Galois groups. This correspondence between ad\`elic
and Galois groups naturally arises within the context of Shimura
varieties. In fact, a Shimura variety is a pro-variety defined
over $\Q$, with a rich ad\`elic group of symmetries. In that
context, the compatibility of the Galois action and the
automorphisms is at the heart of Langlands program. This leads us
to give a reinterpretation of the BC and the $\GL_2$ systems in
the language of Shimura varieties, with the BC system
corresponding to the simplest (zero dimensional) Shimura variety.
In the case of the $\GL_2$ system, we show how the data of
2-dimensional $\Q$-lattices and commensurability can be also
described in terms of elliptic curves together with a pair of
points in the total Tate module, and the system is related to the
Shimura variety of $\GL_2$. This viewpoint suggests considering
our systems as {\em noncommutative pro-varieties} defined over
$\Q$, more specifically as noncommutative Shimura varieties.

\smallskip

We then present our main result, which is the construction of a new
system, whose arithmetic properties fully
incorporate the explicit class field theory for an imaginary quadratic
field $K$, and whose partition function is the Dedekind zeta function of
$K$. The underlying geometric structure is given by commensurability of
1-dimensional $K$-lattices.

\smallskip

This new system can be regarded in two different ways. On the one hand, 
it is a generalization of the BC system of \cite{BC}, when changing the
field from $\Q$ to $K$, and is in fact Morita equivalent to the 
one considered in \cite{LvF}, but with no restriction on the class number. 
On the other hand, it is also a specialization of the
$\GL_2$-system of \cite{CM} to elliptic curves with complex multiplication by
$K$. In this case the ground states can be related to the
non-generic ground states of the $\GL_2$-system, associated to
points $\tau\in \H$ with complex multiplication, and the group of
symmetries is the Galois group of the maximal abelian extension of
$K$.

\smallskip

Here also we show that symmetries by endomorphisms play a crucial
role, as they allow for the action of the class group ${\rm
Cl}(\O)$, so that our results hold for any class field. Since this
complex multiplication (CM) case can be realized as a subgroupoid
of the $\GL_2$-system, it has a natural choice of a rational
subalgebra (an arithmetic structure) inherited from that of the
$\GL_2$-system. This is crucial, in order to obtain the
intertwining of Galois action on the values of extremal states and
action of symmetries of the system.

\smallskip

We summarize and compare the main properties of the three systems
(BC, $\GL_2$, and CM) in the following table.

\bigskip
\bigskip
\begin{center}
\begin{tabular}{|c|c|c|c|}
\hline & & & \\
system & $\GL_1$ & $\GL_2$ & CM \\
\hline & & & \\
partition function & $\zeta(\beta)$ & $\zeta(\beta)\zeta(\beta-1)$
& $\zeta_K(\beta)$ \\
\hline & & & \\
symmetries & $\A^*_f/\Q^*$ & $\GL_2(\A_f)/\Q^*$ & $\A^*_{K,f}/K^*$ \\
\hline & & & \\
automorphisms & $\hat\Z^*$ & $\GL_2(\hat\Z)$ &
$\hat\O^*/\O^*$ \\
\hline & & & \\
endomorphisms & & $\GL_2^+(\Q)$ & ${\rm Cl}(\O)$ \\
\hline & & & \\
Galois group & $\Gal(\Q^{ab}/\Q)$ & $\Aut(F)$ & $\Gal(K^{ab}/K)$ \\
\hline & & & \\
extremal KMS$_\infty$ & $Sh(\GL_1,\pm 1)$ & $Sh(\GL_2,\H^\pm)$ &
$\A^*_{K,f}/K^*$ \\
\hline
\end{tabular}
\end{center}
\bigskip
\bigskip

The paper consists of two parts, with sections
\ref{QScft} and \ref{QlatShi} centered on the relation of the 
BC and $\GL_2$ system to the arithmetic of Shimura varieties, 
and sections \ref{CMsect} and \ref{KMSsect} dedicated to the 
construction of the CM system and its relation to the explicit 
class field theory for imaginary quadratic fields. The two parts 
are closely interrelated, but can also be read independently.

\section{Quantum Statistical
Mechanics and Explicit Class Field Theory}\label{QScft}

  The BC quantum statistical mechanical system \cite{bos-con-CR, BC}
exhibits generators of the maximal abelian extension of $\Q$,
parameterizing ground states (\ie at zero temperature). Moreover, the system
has the remarkable property that these ground states take algebraic
values, when evaluated on a rational subalgebra of the $C^*$-algebra
of observables. The action on these values of the absolute Galois
group factors through the abelianization $\Gal(\Q^{ab}/\Q)$ and is
implemented by the action of the id\`ele class group
as symmetries of the system, via the class field theory isomorphism.
This suggests the intriguing possibility of using the setting of
quantum statistical mechanics to address the problem of explicit class
field theory for other number fields.

\smallskip

In this section we recall some basic notions of quantum statistical
mechanics and of class field theory, which will be used throughout the paper.
We also formulate a general conjectural relation between quantum statistical
mechanics and the explicit class field theory problem for number fields.

\medskip
\subsection*{Quantum Statistical Mechanics}\hfill\medskip

  A quantum statistical mechanical system consists of an algebra of
observables, given by a unital $C^*$-algebra $\cA$, together with a
time evolution, consisting of a 1-parameter group of automorphisms
$\sigma_t$, \((t\in \R)\), whose infinitesimal generator is the
Hamiltonian of the system. The analog of a probability measure,
assigning to every observable a certain average, is given by a state,
namely a continuous linear functional $\varphi: \cA\to \C$
satisfying positivity, $\varphi( x^* x)\geq 0$, for all $x\in \cA$, and
normalization, $\varphi(1)=1$. In the quantum mechanical framework,
the analog of the classical Gibbs measure is given by states
satisfying the KMS condition (\cf \cite{hhw}).

\begin{defn}\label{KMSbetadef}
A triple $(\cA, \sigma_t,\varphi)$ satisfies the
Kubo-Martin-Schwinger (KMS) condition at inverse temperature
$0\leq \beta < \infty$, if, for all $x,y\in \cA$, there exists a
bounded holomorphic function $F_{x,y}(z)$ on the strip $0< {\rm Im}(z)<
\beta$, continuous on the boundary of the strip, such that
\begin{equation}\label{KMSdef}
 F_{x,y}(t)=\varphi(x\sigma_t(y)) \ \ \ \text{ and } \ \ \
F_{x,y}(t+i\beta)=\varphi(\sigma_t(y)x), \ \ \ \ \forall t\in \R.
\end{equation}
\end{defn}

  We also say that $\varphi$ is a KMS$_\beta$ state for \((\cA,
\sigma_t)\). The set $\cK_\beta$ of KMS$_\beta$ states is a compact
convex Choquet simplex \cite[II \S5]{BR} whose set of extreme
points $\E_\beta$ consists of the factor states.

\smallskip

  At $0$ temperature ($\beta = \infty$) the KMS condition
\eqref{KMSdef} says that,
for all $x,y\in \sA$, the function
\begin{equation}\label{Ftab}
F_{x,y}(t)= \varphi(x\, \sigma_t(y))
\end{equation}
extends to a bounded holomorphic function in the upper half plane
$\Hb$. This implies that, in the Hilbert space of the GNS
representation of $\varphi$ (\ie the completion of $\sA$ in the inner
product $\varphi(x^*y)$), the generator $H$ of the one-parameter group
$\sigma_t$ is a positive operator (positive energy condition).
However, this notion of $0$-temperature KMS states is in general too
weak, hence the notion of KMS$_\infty$ states that we shall consider is
the following.

\begin{defn}\label{KMSinftydef}
A state $\varphi$ is a KMS$_\infty$ state for \((\cA,
\sigma_t)\) if it is a weak limit of $\beta$-KMS states for $\beta
\rightarrow \infty$.
\end{defn}

  One can easily see the difference between these two notions
in the case of the trivial evolution
$\sigma_t = {\rm id} \qqq t\in \R$, where any
state has the property that \eqref{Ftab} extends to the upper half
plane (as a constant), while weak limits of
$\beta$-KMS states are automatically tracial states.
With Definition \ref{KMSinftydef} we still obtain a weakly
compact convex set
$\Sigma_\infty$ and we can consider the set $\E_\infty$ of its
extremal points.

\smallskip

  The typical framework for spontaneous symmetry breaking in a
system with a unique phase transition (\cf \cite{haag}) is that the
simplex $\Sigma_\beta$ consists of a single point for $\beta \leqq
\beta_c$ \ie when the temperature is larger than the
critical temperature $T_c$, and is non-trivial (of some higher
dimension in general) when the temperature lowers.
A (compact) group of automorphisms
$G \subset {\rm Aut} ({\mathcal A})$ commuting
with the time evolution,
\begin{equation}
\label{eq5}
\sigma_t \, \alpha_g = \alpha_g \, \sigma_t \qquad \forall \, g \in G
\, , \ t  \in {\mathbb R} \, ,
\end{equation}
is a symmetry group of the system. Such $G$ acts on $\Sigma_{\beta}$
for any $\beta$, hence on the
extreme points ${\mathcal E} (\Sigma_{\beta}) = {\mathcal E}_{\beta}$.
The choice of an equilibrium state $\varphi\in {\mathcal E}_{\beta}$
may break this symmetry to a smaller subgroup given by the isotropy group
$G_\varphi=\{ g \in G \, , \ g \varphi = \varphi \}$.

\smallskip

  The unitary group $\Uc$ of the
fixed point algebra of $\sigma_t$ acts by inner
automorphisms of the dynamical system $({\mathcal A},\sigma_t)$, by
\begin{equation}\label{inneract}
({\rm Ad}u)\,(a):=\,u\,a\,u^\ast \qqq a\in {\mathcal A}\,,
\end{equation}
for all $u\in \Uc$. One can define an action {\em modulo inner} of a
group $G$ on the system $({\mathcal A},\sigma_t)$
as a map $\a:G\to {\rm Aut}({\mathcal A},\sigma_t)$
fulfilling the condition
\begin{equation}\label{innaction}
\a(gh)\,\a(h)^{-1}\,\a(g)^{-1}\in {\rm Inn}({\mathcal A},\sigma_t)
\qqq  g,h \in G\,,
\end{equation}
\ie, as a homomorphism of $G$ to $\Aut(A,\sigma_t)/\Uc$.
The KMS$_\beta$ condition shows that the
{\em inner} automorphisms ${\rm Inn}({\mathcal A},\sigma_t)$
act trivially on KMS$_\beta$ states, hence
\eqref{innaction} induces an action of the group $G$ on the set
$\Sigma_{\beta}$ of KMS$_\beta$ states, for $0<\beta\leq \infty$.

\smallskip

  More generally, one can consider actions {\em by
endomorphisms} (\cf \cite{CM}), where an endomorphism $\rho$ of the
dynamical system $({\mathcal A},\sigma_t)$ is a $\ast$-homomorphism
$\rho:{\mathcal A} \to {\mathcal A}$ commuting with the evolution
$\sigma_t$. There is an induced action of $\rho$ on KMS$_\beta$
states, for $0<\beta< \infty$, given by
\begin{equation}\label{endoKMSact}
\rho^\ast(\varphi):= Z^{-1}\,\varphi \circ \rho \,,\quad Z=\varphi(e),
\end{equation}
provided that $\varphi(e)\neq 0$, where $e=\rho(1)$ is an idempotent
fixed by $\sigma_t$.

\smallskip

  An {\em isometry} $u\in \cA$, $u^\ast\,u=1$, satisfying
$\sigma_t(u)=\lambda^{it} \,u$ for all $t\in \R$ and for some $\lambda
\in\R^\ast_+$, defines an {\em inner} endomorphism ${\rm Ad}u$
of the dynamical system $({\mathcal A},\sigma_t)$, again of the form
\eqref{inneract}. The KMS$_\beta$ condition shows that the
induced action of ${\rm Ad}u$ on $\Sigma_\beta$ is trivial, \cf
\cite{CM}. The induced action
(modulo inner) of a semigroup of endomorphisms of
$(\cA,\sigma_t)$ on the KMS$_\beta$ states in general may
not extend directly to KMS$_\infty$ states (in a nontrivial way), but
it may be defined on $\E_\infty$ by ``warming up and cooling
down'' (\cf \cite{CM}), provided the ``warming up'' map $W_\beta:
\E_\infty \to \E_\beta$ is a bijection between KMS$_\infty$
states (in the sense of Definition \ref{KMSinftydef}) and KMS$_\beta$
states, for sufficiently large $\beta$. The map is given by
\begin{equation}\label{warming}
W_\beta(\varphi)(a)= \frac{
\Tr (\pi_\varphi(a)\,e^{-\beta \,H})}{\Tr(\,e^{-\beta \,H})}
\qqq a\in \cA\,,
\end{equation}
with $H$ the positive energy Hamiltonian, implementing the time
evolution in the representation $\pi_\varphi$ associated to the
extremal KMS$_\infty$ state $\varphi$.

\smallskip

  This type of symmetries, implemented by endomorphisms instead of
automorphisms, plays a crucial role in the theory of superselection
sectors in quantum field theory, developed by Doplicher--Haag--Roberts
(cf.\cite{haag}, Chapter IV).

\smallskip

States on a $C^*$-algebra extend the notion of integration with
respect to a measure in the commutative case. In the case of a
non-unital algebra, the multipliers algebra provides a
compactification, which corresponds to the Stone--\v{C}ech
compactification in the commutative case. A state admits a canonical
extension to the multiplier algebra. Moreover, just as in the
commutative case  one can extend integration to certain classes of
unbounded functions, it is preferable to extend, whenever possible,
the integration provided by a state to certain classes of unbounded
multipliers.

\medskip
\subsection*{Hilbert's 12th problem}\hfill\medskip

  The main theorem of class field theory provides a
classification of finite abelian extensions of a local or global
field $K$ in terms of subgroups of a locally compact abelian group
canonically associated to the field. This is the multiplicative
group $K^*=\GL_1(K)$ in the local non-archimedean case, while in the global case
it is the quotient of the id\`ele class group $C_K$ by the
connected component of the identity. The construction of the group
$C_K$ is at the origin of the theory of id\`eles and ad\`eles.

\smallskip

  Hilbert's 12th problem can be formulated as the question
of providing an explicit description of a set of generators of the
maximal abelian extension $K^{ab}$ of a number field $K$ and of the
action of the Galois group $\Gal(K^{ab}/K)$. This is the maximal
abelian quotient of the absolute Galois group \({\rm Gal}(\bar{K}/K)\)
of \(K\), where \(\bar{K}\) denotes an algebraic closure of \(K\).

\smallskip

  Remarkably, the only cases of number fields for which there is a
complete answer to Hilbert's 12th problem are the construction of
the maximal abelian extension of $\Q$ using torsion points of $\C^*$
(Kronecker--Weber) and the case of imaginary quadratic fields, where
the construction relies on the theory of elliptic curves with complex
multiplication (\cf \eg the survey \cite{St}).

\smallskip

  If $\A_K$ denotes the ad\`eles of a number field $K$ and
$J_K=\GL_1(\A_K)$ is the group of id\`eles of $K$, we write $C_K$ for
the group of id\`ele classes $C_K=J_K/K^*$ and $D_K$ for the
connected component of the identity in $C_K$.

\smallskip

  In the rest of this paper the field \(K\) will be either \(\Q\)
or an imaginary quadratic field \(\Q(\sqrt{-d})\), for some
positive integer \(d >1\). We denote by \(\hat{\Z}\)  the
profinite completion of \(\Z\) and by \(\A_f= \hat{\Z}\otimes\Q\)
the ring of finite adeles of \(\Q\). For any abelian group \(G\),
we denote by \(G_{tors}\) the subgroup of elements of finite
order. For any ring \(R\), we write \(R^*\) for the group of
invertible elements, while \(R^{\times}\) denotes the set of
nonzero elements of \(R\), which is a semigroup if $R$ is an
integral domain. We write \(\mathcal O\) for the ring of algebraic
integers of \(K\). We set \(\hat{\mathcal O}:= (\mathcal O
\otimes\hat{\Z})\) and write $\A_{K,f}=\A_f\otimes_\Q K$ and
\(\J_K =\A_{K,f}^*=\GL_1(\A_{K,f})\). Note that \(K^*\) embeds
diagonally into \(\J_K\). The class field theory isomorphism
provides the canonical identification
\begin{equation}\label{CFTiso}
\theta: \J_K/K^* \xrightarrow{\sim} {\rm Gal}(K^{ab}/K)\,,
\end{equation}
with $K^*$ replaced by $\Q_+^*$ when $K=\Q$.

\medskip
\subsection*{Fabulous states for number fields}\hfill\medskip

  The connection between class field theory and quantum statistical
mechanics can be formulated as the problem of constructing a class
of quantum statistical mechanical systems, whose set of ground states
$\E_\infty$ has special arithmetic properties, because of which we
refer to such states as ``fabulous states''.

\smallskip

  Given a number field $K$, with a choice of an embedding $K\subset
\C$, the ``problem of fabulous states'' consists of constructing
a $C^*$-dynamical system $({\mathcal A}, \sigma_t)$, with an {\em
arithmetic subalgebra} ${\mathcal A}_\Q$ of ${\mathcal A}$,
with the following properties:
\begin{enumerate}
\item The quotient group $G=C_K/D_K$ acts on ${\mathcal A}$ as symmetries
compatible with $\sigma_t$.
\item The states $\varphi \in \sE_\infty$, evaluated on elements of
the arithmetic subalgebra $\cA_\Q$, satisfy:
\begin{itemize}
\item $\varphi(a)\in \overline{K}$, the algebraic
closure of $K$ in $\C$;
\item the
elements of $\{ \varphi(a):\,\, a\in {\mathcal A}_K, \,\, \varphi\in
\sE_\infty\}$ generate $K^{ab}$.
\end{itemize}
\item The class field theory isomorphism
\begin{equation}\label{CFTisoCD}
\theta:C_K/D_K \stackrel{\simeq}{\longrightarrow} \Gal (K^{ab}/K)
\end{equation}
intertwines the actions,
$$ \alpha \circ \varphi = \varphi \circ \theta^{-1}(\alpha), $$
for all $\alpha \in \Gal (K^{ab}/K)$ and for all $\varphi \in
\sE_\infty$.
\end{enumerate}

\smallskip

 In the setting described above
the $C^*$-dynamical system \((\cA, \sigma_t)\) together
with a \(\Q\)-structure compatible with the flow \(\sigma_t\) (\ie
a rational subalgebra $\cA_\Q\subset \cA$ such that
\(\sigma_t(\cA_\Q\otimes\C) = \cA_\Q\otimes\C\)) defines
 a \emph{non-commutative algebraic (pro-)variety} \(X\) over
\(\Q\). The ring $\cA_\Q$ (or \(\cA_\Q\otimes\C\)), which need not
be involutive, is the analog
of the ring of algebraic functions on \(X\) and the set of
extremal KMS\(_{\infty}\)-states is the analog of the set of
points of \(X\). The action of the subgroup of \({\rm
Aut}(\cA,\sigma_ t)\) which takes \(\cA_\Q\otimes\C\) into itself
is analogous to the action of the Galois group on the (algebraic)
values of algebraic functions at points of \(X\).

\smallskip

The analogy illustrated above leads us to speak somewhat loosely
of ``classical points'' of such a noncommutative algebraic
pro-variety. We do not attempt to give a general definition of
classical points, while we simply remark that, for the specific
construction considered here, such a notion is provided by the
zero temperature extremal states.

\medskip

A broader type of question, in a similar spirit, can be formulated
regarding the construction of quantum statistical mechanical systems
with ad\`elic groups of symmetries and the arithmetic properties of its
action on zero temperature extremal KMS states. The case of the
$\GL_2$-system of \cite{CM} fits into this general program.

\section{$\Q$-lattices and noncommutative Shimura varieties}\label{QlatShi}

In this section we recall the main properties of the BC and the
$\GL_2$ system, which will be useful for our main result.

\smallskip

Both cases can be formulated starting with the same geometric notion,
that of commensurability classes of $\Q$-lattices, in dimension one
and two, respectively.

\medskip

\begin{defn}\label{defQlat} A $\Q$-lattice in $\R^n$ is a
 pair $ \, ( \L , \phi) \,$, with $
\L $ a lattice in $\R^n$, and
\begin{equation}\label{phi-label}
\phi :  \Q^n/\Z^n \longrightarrow \Q\L / \L
\end{equation}
a homomorphism of abelian groups.  A $\Q$-lattice is invertible if the map
\eqref{phi-label} is an isomorphism. Two $\Q$-lattices $(\Lambda_1,
\phi_1)$ and $(\Lambda_2, \phi_2)$ are commensurable
if the lattices are commensurable (\ie $\Q\Lambda_1=\Q\Lambda_2$) and the
maps $\phi_1$ and $\phi_2$ agree modulo the sum of the lattices.
\end{defn}

It is essential here that one does not require the homomorphism $\phi$
to be invertible in general.

\smallskip

The set of $\Q$-lattices modulo the equivalence relation of
commensurability and considered up to scaling is best described with the
tools of noncommutative geometry, as explained in \cite{CM}. In fact,
one can first consider the groupoid of the equivalence relation of
commensurability on the set of $\Q$-lattices (not up to scaling). This
is a locally compact \'etale groupoid $\cR$. When considering the quotient
by the scaling action (by $S=\R_+^*$ in the 1-dimensional case, or by
$S=\C^*$ in the 2-dimensional case), the algebra of coordinates
associated to the quotient $\cR/S$ is obtained by restricting the
convolution product of the algebra of $\cR$ to weight zero functions
with $S$-compact support. The algebra obtained this way, which is
unital in the 1-dimensional case, but not in the 2-dimensional case,
has a natural time evolution given by the ratio of the covolumes of a
pair of commensurable lattices.
Every unit $y\in \cR^{(0)}$ of $\cR$ defines a representation
$\pi_y$ by left convolution of the algebra of $\cR$ on the Hilbert
space $\cH_y=\ell^2(\cR_y)$, where $\cR_y$ is the set of elements with
source $y$. This construction passes to the quotient by the scaling
action of $S$. Representations corresponding to
points that acquire a nontrivial automorphism group will no longer be
irreducible. If the unit $y\in \cR^{(0)}$
corresponds to an invertible $\Q$-lattice, then $\pi_y$ is a positive
energy representation.

\smallskip

In both the 1-dimensional and the 2-dimensional case, the set of
extremal KMS states at low temperature is given by a classical
ad\`elic quotient, namely, by the Shimura varieties for $\GL_1$ and
$\GL_2$, respectively, hence we argue here that the noncommutative
space describing commensurability classes of $\Q$-lattices up to scale
can be thought of as a {\em
noncommutative Shimura variety}, whose set of classical points is the
corresponding classical Shimura variety.

\smallskip

In both cases, a crucial step for the arithmetic properties of the
action of symmetries on extremal KMS states at zero temperature is
the choice of an arithmetic subalgebra
of the system, on which the ground states are evaluated. Such choice
gives the underlying noncommutative space a more rigid structure, of
``noncommutative arithmetic variety''.

\medskip
\subsection*{Tower Power}\hfill\medskip

If \(V\) is an algebraic variety -- or a scheme or a stack -- over
a field \(k\), a ``tower'' \(\cT\) over \(V\) is a family \(V_i \)
\((i \in \I)\) of finite (possibly branched) covers of \(V\) such
that for any \(i,j\in \I\), there is a \(l\in \I\) with \(V_l\) a
cover of \(V_i\) and \(V_j\).  Thus,  \(\I\) is a partially ordered
set. In case of a
tower over a pointed variety \((V,v)\), one fixes a point \(v_i\)
over \(v\) in each \(V_i\). Even though \(V_i\) may not be
irreducible, we shall allow
ourselves to loosely refer to \(V_i\) as a variety.
It is convenient to view a ``tower'' \(\cT\) as a category
\(\mathcal C\) with objects
\((V_i \to V)\) and morphisms \({\rm Hom}(V_i,V_j)\) being maps of covers of
\(V\). One has the group \({\rm Aut}_{\cT}(V_i)\)  of
invertible self-maps of \(V_i\) over \(V\) (the group of deck
transformations); the deck transformations
 are not required to preserve the points \(v_i\). There is a (profinite)
 group of symmetries associated to a tower, namely
\begin{equation}\label{proG}
 \cG:= {\varprojlim}_{i} \text{Aut}_{\cT}(V_i).
\end{equation}

\smallskip

The simplest example of a tower is the ``fundamental group'' tower
associated with a (smooth connected) complex algebraic variety
\((V,v)\)  and its
universal covering \((\tilde{V}, \tilde{v})\).  Let \(\mathcal C\) be
the category
of all finite \'etale (unbranched) intermediate covers \(\tilde{V} \to
W \to V\) of \(V\) . In this case, the symmetry group \(\cG\)
of \eqref{proG} is the algebraic fundamental group of \(V\); it is
also the profinite completion of the (topological) fundamental group
\(\pi_1(V, v)\). Simple variants of this example
include allowing controlled ramification. Other examples of towers are
those defined by iteration of self maps of algebraic varieties.

\smallskip

For us, the most important examples of ``towers'' will be the
cyclotomic tower and the modular tower. Another very interesting case of
towers is that of more general Shimura varieties. These, however, will
not be treated in this paper. (For a more general treatment of 
noncommutative Shimura varieties see \cite{HaPau}.)

\medskip
\subsection*{The cyclotomic tower and the BC system}\hfill\medskip

In the case of \(\Q\), an explicit description of \(\Q^{ab}\) is
provided by the Kronecker--Weber theorem.
This shows that the field
\(\Q^{ab}\) is equal to
  \(\Q^{cyc}\), the field obtained by attaching all roots of unity to
  \(\Q\). Namely, \(\Q^{ab}\)  is obtained by attaching the values of
  the exponential function \(exp(2\pi i z)\) at the torsion
  points of the circle group \(\R/{\Z}\). Using the
  isomorphism of abelian groups
  \(\bar{\Q}^*_{tors} \cong \Q/{\Z}\)
and the identification \({\rm Aut}(\Q/{\Z}) = {\rm
    GL}_1(\hat{\Z}) = \hat{\Z}^*\), the restriction to
  \(\bar{\Q}^*_{tors}\) of the natural
  action of \({\rm Gal}(\bar{\Q}/{\Q})\) on \(\bar{\Q}^*\) factors as
\[
{\rm Gal}(\bar{\Q}/{\Q}) \to {\rm Gal}(\bar{\Q}/{\Q})^{ab} = {\rm
  Gal}({\Q}^{ab}/{\Q}) \xrightarrow{\sim} \hat{\Z}^*.
\]

\smallskip

Geometrically, the above setting can be understood in terms of the
{\em cyclotomic tower}. This has
base Spec~\(\Z = V_1\). The family is Spec~\(\Z[\zeta_n] = V_n\)
where \(\zeta_n\) is a primitive \(n\)-th root of unity (\(n \in
\N^*\)). The set \({\rm Hom}~(V_m \to V_n)\), non-trivial for \(n|m\),
 corresponds to the map \(\Z[\zeta_n]
\hookrightarrow \Z[\zeta_m]\) of rings. The group \({\rm Aut}
(V_n) = {\rm GL}_1(\Z/{n\Z})\) is the Galois group \({\rm Gal}
(\Q(\zeta_n)/{\Q})\). The group of symmetries \eqref{proG} of the
tower is then
\begin{equation}\label{hatZgroup}
 \cG = \varprojlim_n \text{GL}_1(\Z/{n\Z}) =
\text{GL}_1(\hat{\Z}),
\end{equation}
which is isomorphic to the Galois group \({\rm Gal}(\Q^{ab}/\Q)\) of
the maximal abelian extension of \(\Q\).

\smallskip

The classical object that we consider, associated to the cyclotomic
tower, is the Shimura variety given by the ad\`elic quotient
\begin{equation}\label{ShGL1}
Sh(\GL_1,\{\pm 1\}) = \GL_1(\Q) \backslash (\GL_1(\A_f)\times \{
\pm 1 \}) = \A_f^*/\Q^*_+.
\end{equation}

\smallskip

Now we consider the space of 1-dimensional $\Q$-lattices up to scaling
modulo commensurability. This can be described as follows
(\cite{CM}).

\smallskip

In one dimension, every $\Qb$-lattice is of the form
\begin{equation}\label{1dQlatrho}
 ( \L , \phi) \, = (\lambda\, \Z,\lambda\,\rho),
\end{equation}
for some $\lambda>0$ and some
$\rho \in \Hom(\Q/\Z,\Q/\Z)$. Since we can identify
$\Hom(\Q/\Z,\Q/\Z)$
endowed with the topology of pointwise convergence with
\begin{equation}\label{RhatZ}
\Hom (\Q/\Z,\Q/\Z) = \varprojlim_n \Z/n\Z = \hat \Z\,,
\end{equation}
we obtain that the algebra $C(\hat \Z)$ is the algebra of coordinates
of the space of 1-dimensional $\Q$-lattices up to scaling.
The group $\hat\Z$ is  the Pontrjagin dual of $\Q/\Z$, hence we also
have an identification $C(\hat \Z)=C^*(\Q/\Z)$.

\smallskip

The group of deck
transformations $\cG=\hat\Z^*$ of the cyclotomic tower acts by
automorphisms on the algebra of coordinates $C(\hat \Z)$. In addition
to this action, there is a semigroup action of $\N^\times=\Z_{>0}$
implementing the commensurability relation. This is given by
endomorphisms that move vertically across the levels of the cyclotomic
tower. They are given by
\begin{equation}\label{alphax-cycl}
\alpha_n(f) (\rho)= f(n^{-1} \rho), \ \ \ \forall \rho \in
n\hat\Z.
\end{equation}
Namely, $\alpha_n$ is the isomorphism of $C(\hat\Z)$ with the
reduced algebra $C(\hat\Z)_{\pi_n}$ by the projection $\pi_n$
given by the characteristic function of $n\hat\Z\subset \hat\Z$.
Notice that the action \eqref{alphax-cycl} cannot be restricted to
the set of invertible $\Q$-lattices, since the range of $\pi_n$ is
disjoint from them.

\smallskip

The algebra of coordinates $\cA_1$ on the {\em noncommutative}
space of equivalence classes of 1-dimensional $\Q$-lattices modulo
scaling, with respect to the equivalence relation of
commensurability, is given then by the semigroup crossed product
\begin{equation}\label{BCalg}
 \cA=C(\hat\Z)\rtimes_\alpha \N^\times.
\end{equation}

\smallskip

Equivalently, we are considering the convolution algebra of the
groupoid $\cR_1/\R^*_+$ given by the quotient by scaling of the
groupoid of the equivalence relation of commensurability on
1-dimensional $\Q$-lattices, namely, $\cR_1/\R^*_+$ has as algebra of
coordinates the functions
$f(r,\rho)$, for $\rho\in \hat\Z$ and $r\in \Q^*$ such that $r\rho\in
\hat\Z$, with the convolution product
\begin{equation}\label{convol1}
 f_1 * f_2 \, (r,\rho)= \sum f_1(rs^{-1}, s\rho) f_2(s,\rho),
\end{equation}
and the adjoint $f^*(r,\rho)=\overline{f(r^{-1},r\rho)}$.

\smallskip

This is the $C^*$-algebra of the Bost--Connes (BC) system
\cite{BC}. It was originally defined as a Hecke
algebra for the almost normal pair of solvable groups
$P_\Z^+\subset P_\Q^+$, where $P$ is the algebraic $ax+b$ group and
$P^+$ is the restriction to $a>0$ (\cf
\cite{BC}). It has a natural time evolution $\sigma_t$ determined by
the regular representation of this Hecke algebra, which is of type
III$_1$. The time evolution depends upon the ratio of the lengths of
$P_\Z^+$ orbits on the left and right cosets.

\smallskip

As a set, the space of commensurability
classes of 1-dimensional $\Q$-lattices up to scaling can also be
described by the quotient
\begin{equation}\label{quotient1NC}
\GL_1(\Q)\backslash \A^\cdot / \R^*_+ = \GL_1(\Q)\backslash (\A_f
\times \{ \pm 1 \}),
\end{equation}
where $\A^\cdot:=\A_f\times \R^*$ is the set of ad\`eles with nonzero
archimedean component. Rather than considering this quotient set
theoretically, we regard it as a noncommutative space, so as to be
able to extend to it the ordinary tools of geometry that can be
applied to the ``good'' quotient \eqref{ShGL1}.

\smallskip

The noncommutative algebra of coordinates of \eqref{quotient1NC}
is the crossed product
\begin{equation}\label{crossQA1}
C_0(\A_f) \rtimes \Q^*_+ .
\end{equation}
This is Morita equivalent to the algebra \eqref{BCalg}. In fact,
\eqref{BCalg} is obtained as a full corner of \eqref{crossQA1},
$$ C(\hat\Z)\rtimes \N^\times = \left( C_0(\A_f) \rtimes \Q^*_+
\right)_\pi, $$
by compression with the projection $\pi$ given by the characteristic
function of $\hat\Z \subset \A_f$ (\cf \cite{laca-end}).

\smallskip

The quotient \eqref{quotient1NC} with its noncommutative algebra of
coordinates \eqref{crossQA1} can then be thought of as the {\em
noncommutative Shimura variety}
\begin{equation}\label{Sh1nc}
Sh^{(nc)}(\GL_1,\{ \pm 1\}):= \GL_1(\Q)\backslash (\A_f \times
\{\pm 1 \})= \GL_1(\Q)\backslash \A^\cdot /\R^*_+,
\end{equation}
whose set of classical points is the well behaved quotient
\eqref{ShGL1}.

\smallskip

This has a {\em compactification}, obtained by replacing
$\A^\cdot$ by $\A$, as in \cite{Co-zeta},
\begin{equation}\label{compSh1nc}
 \overline{Sh^{(nc)}}(\GL_1,\{ \pm 1 \}) = \GL_1(\Q)\backslash \A
/\R^*_+.
\end{equation}
The compactification consists of adding the trivial lattice (with
a possibly nontrivial $\Q$-structure).

\smallskip

One can also consider the noncommutative space dual to
\eqref{compSh1nc}, under the
duality given by taking the crossed
product by the time evolution. This
is the noncommutative space that gives the spectral realization of
the zeros of the Riemann zeta function in \cite{Co-zeta}. It is a
principal $\R^*_+$-bundle over the noncommutative space
\begin{equation}\label{compact1NC}
\GL_1(\Q)\backslash \A/ \R^*_+.
\end{equation}

\medskip
\subsection*{Arithmetic structure of the BC system}\hfill\medskip

The results of \cite{BC} show that the Galois theory of the cyclotomic
field $\Q^{cycl}$ appears naturally in the BC system when considering
the action of the group of symmetries of the system on the extremal
KMS states at zero temperature.

\smallskip

In the case of 1-dimensional $\Q$-lattices up to scaling, the
algebra of coordinates $C(\hat \Z)$ can be regarded as the
algebra of {\em homogeneous functions of weight zero} on the space of
1-dimensional $\Q$-lattices. As such, it has a natural choice of an
arithmetic subalgebra.

\smallskip

This is obtained in \cite{CM} by considering functions on the
space of 1-dimensional $\Q$-lattices of the form
\begin{equation}\label{epsilonk}
\epsilon_{1,a} ( \L , \phi)= \sum_{y\in  \Lambda +\phi(a)} y^{-1},
\end{equation}
for any $a\in \Q/\Z$. This is well defined, for $\phi(a)\neq 0$,
using the summation $\lim_{N \rightarrow \infty}\sum_{-N}^N$.
numbers. One can then form the weight zero functions
\begin{equation}\label{eka}
e_{1,a}:= c \, \epsilon_{1,a},
\end{equation}
where $c(\Lambda)$ is proportional to the covolume $|\Lambda|$ and
normalized so that $(2\pi \sqrt{-1}) c(\Z)=1$. One considers the
rational subalgebra $\cA_{1,\Q}$ of \eqref{BCalg} generated by the
functions $e_{1,a}(r,\rho):=e_{1,a}(\rho)$ and by the 
functions $\mu_n(r,\rho)=1$ for $r=n$ and zero
otherwise, that implement the semigroup action of $\N^\times$ in 
\eqref{BCalg}.

\smallskip

As proved in \cite{CM}, the algebra $\cA_{1,\Q}$ is the same as the
rational subalgebra considered in \cite{BC}, generated over $\Q$ by
the $\mu_n$ and the exponential functions
\begin{equation}\label{errho}
e(r)( \rho):=\exp(2\pi i \rho(r)), \ \ \  \text{ for } \rho
\in \Hom(\Q/\Z,\Q/\Z), \ \ \text{ and  } r\in \Q/\Z,
\end{equation}
with relations $e(r+s)=e(r)e(s)$, $e(0)=1$, $e(r)^*=e(-r)$,
$\mu_n^*\mu_n =1$, $ \mu_k \mu_n=\mu_{kn} $,
and
\begin{equation}\label{muner}
 \mu_n e(r) \mu_n^* = \frac{1}{n} \sum_{ns=r} e(s).
\end{equation}
The $C^*$-completion of $\cA_{1,\Q}\otimes\C$ gives \eqref{BCalg}.

\smallskip

The algebra \eqref{BCalg} has irreducible representations on the
Hilbert space $\cH=\ell^2(\N^\times)$, parameterized by
elements $\alpha\in \hat\Z^*=\GL_1(\hat\Z)$. Any such element
defines an embedding $\alpha: \Q^{cycl}\hookrightarrow \C$ and the
corresponding representation is of the form
\begin{equation}\label{repGL1}
\pi_\alpha (e(r))\, \epsilon_k = \alpha(\zeta_r^k) \, \epsilon_k
\ \ \ \ \pi_\alpha (\mu_n)\,  \epsilon_k =  \epsilon_{nk}.
\end{equation}
The Hamiltonian implementing the time evolution $\sigma_t$ on
$\cH$ is of the form $H\,\epsilon_k=\log k \,\,\, \epsilon_k$
and the partition function of the BC system is then the
Riemann zeta function
$$ Z(\beta)=\Tr\left( e^{-\beta H} \right) = \sum_{k=1}^\infty
k^{-\beta} = \zeta(\beta.).$$

\smallskip

The set $\E_\beta$ of extremal KMS-states of the BC system
enjoys the following properties (\cf \cite{BC}):

\begin{itemize}
\item $\E_\beta=\cK_\beta$ is a singleton for all \(0 <\beta
\le 1\). This unique KMS state takes values
\[ \varphi_{\beta}(e(m/n)) = f_{-\beta+1}(n)/f_1(n), \]
where
\[ f_k(n)=\sum_{d|n} \mu(d) (n/d)^k, \]
with $\mu$ the M\"obius function, and $f_1$ is the Euler totient
function.
\item  For $1<\beta\leq \infty$, elements of $\E_\beta$ are indexed by
the classes of invertible $\Q$-lattices $\rho\in
\hat\Z^*=\GL_1(\hat\Z)$, hence by the classical points \eqref{ShGL1}
of the noncommutative Shimura variety \eqref{Sh1nc},
\begin{equation}\label{quotient1C}
\E_\beta \cong \GL_1(\Q)\backslash \GL_1(\A) / \R^*_+ \cong
C_\Q/D_\Q \cong \J_\Q/\Q^*_+.
\end{equation}
In this range of temperatures, the values of states
$\varphi_{\beta,\rho} \in \E_\beta$ on the elements $e(r)\in \cA_{1,\Q}$
is given, for $1< \beta < \infty$ by polylogarithms evaluated at roots
of unity, normalized by the Riemann zeta function,
\[
\varphi_{\beta,\rho}(e(r)) =
\frac{1}{\zeta(\beta)}\sum_{n=1}^{\infty}n^{-\beta}
\rho (\zeta_r^{k}).
\]
\item The group $\GL_1(\hat\Z)$ acts by automorphisms of the
system. The induced action of  \(\text{GL}_1(\hat{\Z})\)
on the set of extreme KMS states below critical temperature is free
and transitive.
\item  The vacuum states \((\beta = \infty)\) are \emph{fabulous
states} for the field $K=\Q$, namely $\varphi(\cA_{1,\Q})\subset
\Q^{cycl}$ and the class
field theory isomorphism $\theta: \Gal(\Q^{cycl}/\Q)
\stackrel{\cong}{\to} \hat\Z^*$ intertwines the Galois action on
values with the action of $\hat\Z^*$ by symmetries,
\begin{equation}\label{GaloisEq1}
 \gamma \, \varphi(x) = \varphi( \theta(\gamma)\, x),
\end{equation}
for all $\varphi\in \cE_\infty$, for all $\gamma\in
\Gal(\Q^{cycl}/\Q)$ and for all $x\in \cA_{1,\Q}$.
\end{itemize}

\medskip
\subsection*{The modular tower and the $\GL_2$-system}\hfill\medskip

Modular curves arise as moduli spaces of elliptic curves endowed
with additional level structure. Every congruence subgroup
\(\Gamma'\) of \(\Gamma={\rm
  SL}_2(\Z)\) defines a modular curve \(Y_{\Gamma'}\);
we denote by \(X_{\Gamma'}\) the smooth compactification of the
affine curve \(Y_{\Gamma'}\) obtained by adding cusp points.
Especially important among these are the
 modular curves \(Y(n)\) and \(X(n)\) corresponding to the
principal congruence subgroups \(\Gamma(n)\) for \(n \in \N^*\). Any
\(X_{\Gamma'}\) is dominated by an \(X(n)\). We refer to \cite{kula,
  Sh} for more details. We have the following descriptions of the
modular tower.

\smallskip

{\emph{Compact version:}}
The base is \(V=\P^1\) over \(\Q\). The family is given by the modular
curves  \(X(n)\), considered over the cyclotomic field \(\Q(\zeta_n)\)
\cite{Mi}. We note that \({\rm GL}_2(\Z/{n\Z})/{\pm 1}\) is
the group of automorphisms of the projection
\(V_n = X(n) \to X(1) = V_1 =V\).
Thus, we have
\begin{equation}\label{GL2Z}
\cG = \text{GL}_2 (\hat{\Z})/{\pm 1}= \varprojlim_n \text{GL}_2
(\Z/n\Z) / \{ \pm 1 \}.
\end{equation}
{\emph{Non-compact version:}}
The open modular curves \(Y(n)\) form a tower with base the
\(j\)-line \({\rm Spec}~\Q[j] =\A^1 =V_1 -\{\infty\}\).
The ring of modular functions is the union of the rings of functions
of the \(Y(n)\), with coefficients in \(\Q(\zeta_n)\) \cite{kula}.

\smallskip

This shows how the modular tower is a natural geometric way of
passing from $\GL_1(\hat\Z)$ to $\GL_2(\hat \Z)$. The
formulation that is most convenient in our setting is the one given in
terms of Shimura varieties. In fact, rather than the
modular tower defined by the projective limit
\begin{equation}\label{projYn}
Y= \varprojlim_n Y(n)
\end{equation}
of the modular curves $Y(n)$, it is better for
our purposes to consider the Shimura variety
\begin{equation}\label{ShGL2}
Sh(\H^\pm,\GL_2) = \GL_2(\Q) \backslash (\GL_2(\A_f)\times \H^{\pm}) =
\GL_2(\Q) \backslash \GL_2(\A) /\C^*,
\end{equation}
of which \eqref{projYn} is a connected component. In fact, it is
well known that, for arithmetic purposes, it is always better to
work with nonconnected rather than with connected Shimura
varieties (\cf \eg \cite{Mi}). The simple reason why it is
necessary to pass to the nonconnected case is the following. The
varieties in the tower are arithmetic varieties defined over
number fields. However, the number field typically changes along
the levels of the tower ($Y(n)$ is defined over the cyclotomic
field $\Q(\zeta_n)$). Passing to nonconnected Shimura varieties
allows precisely for the definition of a canonical model where the
whole tower is defined over the same number field.

\smallskip

This distinction is important to our viewpoint, since we want to work
with noncommutative spaces endowed with an arithmetic structure,
specified by the choice of an arithmetic subalgebra.

\smallskip

Every 2-dimensional $\Q$-lattice can be described by data
\begin{equation}\label{2dQlatrho}
(\Lambda,\phi) = (\lambda(\Z+\Z z),\lambda\alpha),
\end{equation}
for some $\lambda\in\C^*$, some $z\in \H$, and $\alpha\in
M_2(\hat\Z)$ (using the basis $(1, -z)$ of $\Z+\Z z$
as in (87) \cite{CM} to view $\alpha$ as a map $\phi$).
The diagonal action of $\Gamma= \SL_2(\Z)$
yields isomorphic $\Q$-lattices, and (\cf (87) \cite{CM})
the space of 2-dimensional $\Q$-lattice up to
scaling can be identified with the quotient
\begin{equation}\label{2Qlatscal}
\Gamma\backslash(\M_2(\hat\Z)\times \H).
\end{equation}
The relation of commensurability is implemented by the partially
defined action of $\GL_2^+(\Q)$ on \eqref{2Qlatscal}.

\smallskip

The groupoid $\cR_2$ of the
commensurability relation on 2-dimensional $\Q$-lattices not up
to scaling (\ie the dual space) has as algebra of coordinates
the convolution algebra of $\Gamma\times \Gamma$-invariant functions on
\begin{equation}\label{tildeU}
\tilde\cU =\{ (g,\alpha,u)\in \GL_2^+(\Q)\times M_2(\hat\Z)\times
\GL_2^+(\R)\, | \, g\alpha\in M_2(\hat\Z) \}.
\end{equation}
Up to Morita equivalence, this can also be described as the crossed
product
\begin{equation}\label{tildeAlgShGL2NC}
 C_0(M_2(\A_f) \times \GL_2(\R)) \rtimes \GL_2(\Q).
\end{equation}

\smallskip

When we pass to $\Q$-lattices up to scaling, we take the quotient
$\cR_2/\C^*$.

\smallskip

If $(\Lambda_k,\phi_k)$ $k=1,2$ are a pair of commensurable
2-dimensional $\Q$-lattices, then for any $\lambda\in \C^*$, the
$\Q$-lattices $(\lambda \Lambda_k, \lambda \phi_k)$ are also
commensurable, with
$$ r(g,\alpha, u\lambda)=\lambda^{-1} r(g,\alpha,u). $$
However, the action of $\C^*$ on $\Q$-lattices is not free due to
the presence of lattices $L=(0,z)$, where $z\in \Gamma\backslash \H$
has nontrivial automorphisms.

\smallskip

Thus, the quotient $Z=\cR_2/\C^*$ is no longer a groupoid. This
can be seen in the following simple example. Consider the two
$\Q$-lattices $(\alpha_1,z_1)=(0,2i)$ and $(\alpha_2,z_2)=(0,i)$.
The composite
$((\alpha_1,z_1),(\alpha_2,z_2))\circ((\alpha_2,z_2),(\alpha_1,z_1))$
is equal to the identity $((\alpha_1,z_1),(\alpha_1,z_1))$. We can
also consider the composition
$(i(\alpha_1,z_1),i(\alpha_2,z_2))\circ((\alpha_2,z_2),(\alpha_1,z_1))$,
where $i(\alpha_2,z_2)=(\alpha_2,z_2)$, but this is not the
identity, since $i(\alpha_1,z_1)\neq (\alpha_1,z_1)$.

\smallskip

However, we can still consider the convolution algebra of $Z$, by
restricting the convolution product of $\cR_2$ to homogeneous
functions of weight zero with $\C^*$-compact support, where a
function $f$ has weight $k$ if it satisfies
$$ f(g,\alpha,u\lambda)=\lambda^k f(g,\alpha,u), \ \ \
\forall \lambda \in \C^*. $$

\smallskip

This is the analog of the description \eqref{BCalg} for the
1-dimensional case. The noncommutative algebra of coordinates $\cA_2$
is thus given by a Hecke algebra of functions  on
\begin{equation}\label{Uset}
\cU=\{ (g,\alpha,z)\in \GL_2^+(\Q)\times M_2(\hat\Z)\times \H,
\,\, g\alpha\in M_2(\hat\Z)\}
\end{equation}
invariant under the $\Gamma\times \Gamma$ action
\begin{equation}\label{Gamma2act}
 (g,\alpha, z)\mapsto (\gamma_1 g \gamma_2^{-1}, \gamma_2 \alpha,
\gamma_2(z)),
\end{equation}
with convolution
\begin{equation}\label{Hecke-convol}
(f_1* f_2)\, (g,\alpha,z) = \sum_{s\in \Gamma\backslash \GL_2^+(\Q),\,
s\alpha\in M_2(\hat\Z)}
f_1(gs^{-1},s\alpha,s(z)) \, f_2(s,\alpha,z)
\end{equation}
and adjoint $f^*(g,\alpha,z)=\overline{f(g^{-1},g\alpha,g(z))}$. This
contains the classical Hecke operators
(\cf (128) \cite{CM}). The time evolution
determined by the ratio of covolumes of pairs of commensurable
$\Q$-lattices is given by
\begin{equation}\label{GL2sigmat}
 \sigma_t(f)(g,\alpha,\tau)= \det(g)^{it}
f(g,\alpha,\tau),
\end{equation}
where, for the pair of commensurable $\Q$-lattices
associated to  $(g,\alpha,\tau)$, one has
\begin{equation}\label{GL2cov}
\det(g)= {\rm covolume}(\Lambda')/{\rm covolume}(\Lambda).
\end{equation}

\medskip

We now give a description closer to \eqref{crossQA1}, which shows
that again we can interpret the space of commensurability classes
of 2-dimensional $\Q$-lattices up to scaling as a noncommutative
version of the Shimura variety \eqref{ShGL2}. More precisely, we
give a reinterpretation of the notion of 2-dimensional
$\Q$-lattices and commensurability, which may be useful in
generalising our work to other Shimura varieties.

\smallskip

Implicit in what follows is an isomorphism between \(\Q/{\Z}\) and
the roots of unity \(\mu(\C)\) in \(\C\); for instance, this could
be given by the exponential function \(e^{2\pi i z}\).

\begin{prop}\label{QlatTate}
The data of a 2-dimensional $\Q$-lattice up to scaling are
equivalent to the data of an elliptic curve $E$, together with a
pair of points $\xi=(\xi_1,\xi_2)$ in the cohomology $H^1(E,
\hat{\Z})$. Commensurability of 2-dimensional $\Q$-lattices up to
scale is then implemented by an isogeny of the corresponding
elliptic curves, with the elements $\xi$ and $\xi'$ related via
the induced map in cohomology.
\end{prop}

\proof The subgroup \(\Q\Lambda/{\Lambda}\) of \(\C/{\Lambda}=E\)
is the torsion subgroup \(E_{tor}\) of the elliptic curve \(E\).
Thus, one can rewrite the map \(\phi\)  as a map \(\Q^2/{\Z^2} \to
E_{tor}\). Using the canonical isomorphism \(E[n]
\xrightarrow{\sim}H^1(E, \Z/{n\Z})\), for $E[n]=\Lambda/n\Lambda$
the $n$-torsion points of $E$, one can interpret \(\phi\) as a map
\(\Q^2/{\Z^2} \to H^1(E, \Q/{\Z})\).

\smallskip

By taking \(\Hom (\Q/{\Z}, -)\), the map \(\phi\) corresponds to a
\(\hat{\Z}\)-linear map
\begin{equation}\label{phiH1}
\hat{\Z}\oplus\hat{\Z} \to H^1(E, \hat{\Z}),
\end{equation}
or to a choice of two elements of the latter. In fact, we use here
the identification \( H^1(E, \hat{\Z}) \cong TE \), where $TE$ is
the total Tate module
\begin{equation}\label{torEhatZ}
\Lambda \otimes \hat\Z = \varprojlim_n E[n]=TE,
\end{equation}
so that \eqref{phiH1} gives a cohomological formulation of the
$\hat\Z$-linear map $\phi: \hat{\Z}\oplus\hat{\Z} \to \Lambda
\otimes \hat\Z$. Commensurability of $\Q$-lattices up to scale is
rephrased as the condition that the elliptic curves are isogenous
and the points in the Tate module are related via the induced map
in cohomology.

\endproof

Another reformulation uses the Pontrjagin duality between
profinite abelian groups and discrete torsion abelian groups given
by \(\Hom(-, \Q/{\Z})\). This reformulates the datum \(\phi\) of a
$\Q$-lattice as a \(\hat{\Z}\)-linear map \(
\Hom(\Q\Lambda/\Lambda,\Q/\Z) \to \hat{\Z}\oplus\hat{\Z}\), which
is identified with \(\Lambda\otimes\hat{\Z} \to
\hat{\Z}\oplus\hat{\Z}\). Here we use the fact that \(\Lambda\)
and \(\Lambda\otimes\hat{\Z} \cong H^1(E, \hat{\Z})\) are both
self-dual (Poincar\'e duality of \(E\)). In this dual formulation
commensurability means that the two maps agree on the intersection
of the two commensurable lattices,
$(\Lambda_1\cap\Lambda_2)\otimes\hat\Z$.

\medskip

With the formulation of Proposition \ref{QlatTate}, can then give
a new interpretation of the result of Proposition 43 of \cite{CM},
which shows that the space of commensurability classes of
2-dimensional $\Q$-lattices up to scaling is described by the
quotient
\begin{equation}\label{ShGL2NC}
Sh^{(nc)}(\H^{\pm},\GL_2):=\GL_2(\Q) \backslash (M_2(\A_f)\times
\H^{\pm}).
\end{equation}
In fact, the data $(\Lambda,\phi)$ of a $\Q$-lattice in $\C$ are
equivalent to data $(E,\eta)$ of an elliptic curve $E=\C/\Lambda$
and an $\A_f$-homomorphism
\begin{equation}\label{Afhom}
\eta : \Q^2 \otimes \A_f \to \Lambda \otimes \A_f,
\end{equation}
with $\Lambda\otimes \A_f = (\Lambda \otimes \hat\Z)\otimes \Q$,
where we can identify \(\Lambda\otimes\hat{\Z}\) with the total
Tate module of $E$, as in \eqref{torEhatZ}. Since the $\Q$-lattice
need not be invertible, we do not require that $\eta$ be an
$\A_f$-isomorphism (\cf \cite{Mi}).

\smallskip

The commensurability relation between $\Q$-lattices corresponds to
the equivalence $(E,\eta)\sim (E',\eta')$ given by an isogeny $g:
E \to E'$ and $\eta'=(g\otimes 1)\circ\eta$. Namely, the
equivalence classes can be identified with the quotient of
$M_2(\A_f)\times \H^{\pm}$ by the action of $\GL_2(\Q)$,
$(\rho,z)\mapsto (g\rho, g(z))$.

\smallskip

Thus, \eqref{ShGL2NC} describes a noncommutative Shimura variety
which has the Shimura variety
\eqref{ShGL2} as the set of its classical points. The results of
\cite{CM} show that, as in the case of the BC system, one recovers the
classical points from the low temperature extremal KMS states. We
shall return to this in the next section.

\smallskip

In this case, the
``compactification'', analogous to passing from
\eqref{quotient1NC} to
\eqref{compact1NC}, corresponds to the replacing
\eqref{ShGL2NC} by the noncommutative space
\begin{equation}\label{compShGL2NC}
\overline{Sh^{(nc)}}(\H^{\pm},\GL_2):=\GL_2(\Q) \backslash (M_2(\A_f)\times
\P^1(\C))\sim\GL_2(\Q) \backslash M_2(\A)/\C^*.
\end{equation}
This corresponds to allowing degenerations of the underlying lattice
in $\C$ to a pseudolattice (\cf \cite{Man1}), while maintaining the
$\Q$-structure (\cf \cite{CM}). The ``invertible part''
\begin{equation}\label{inv-compShGL2NC}
\GL_2(\Q) \backslash (\GL_2(\A_f)\times \P^1(\R))
\end{equation}
of the ``boundary'' gives the noncommutative modular tower considered
in \cite{CDS}, \cite{MM}, and \cite{Pau}, so that the full space
\eqref{compShGL2NC} appears to be the most natural candidate for
the geometry underlying the construction of a
quantum statistical mechanical system adapted to the case of both
imaginary and real quadratic fields (\cf \cite{Man1} \cite{Man2}).

\smallskip

In the case of the classical Shimura varieties, the relation between
\eqref{ShGL2} and \eqref{ShGL1} is given by ``passing to components'',
namely we have (\cf \cite{Mi})
\begin{equation}\label{components}
\pi_0(Sh(\GL_2,\H^\pm))=Sh(\GL_1,\{\pm 1\}).
\end{equation}
In fact, the operation of taking connected components of \eqref{ShGL2}
is realized by the map
\begin{equation}\label{pi0Sh}
{\rm sign} \times \det: Sh(\GL_2,\H^{\pm}) \to
Sh(\GL_1,\{\pm 1\}).
\end{equation}

\medskip
\subsection*{Arithmetic properties of the
$\GL_2$-system}\hfill\medskip

The result proved in \cite{CM} for the $\GL_2$-system shows
that the action of symmetries on the extremal
KMS states at zero temperature is now related to the Galois theory of
the field of modular functions.

\smallskip

Since the arithmetic subalgebra for the BC system was obtained by
considering weight zero lattice functions of the form \eqref{eka}, it
is natural to expect that the analog for the $\GL_2$-system will
involve lattice functions given by the Eisenstein series, suitably
normalized to weight zero, according to the analogy developed by
Kronecker between trigonometric and elliptic functions, as outlined by
A.Weil in \cite{WeilEll}. This suggests that modular functions should
appear naturally in the arithmetic subalgebra $\cA_{2,\Q}$ of the
$\GL_2$-system, but that requires working with unbounded multipliers.

\smallskip

This is indeed the case for the arithmetic subalgebra $\cA_{2,\Q}$
defined in \cite{CM}, which we now recall.

\smallskip

Let $F$ be the modular field, namely the field of modular functions
over $\Q^{ab}$ (\cf \eg \cite{Lang}). This is the union of the fields
$F_N$ of modular functions of level $N$ rational over the cyclotomic
field $\Q(\zeta_N)$, that is, such that the $q$-expansion at a cusp has
coefficients in the cyclotomic field $\Q(\zeta_N)$.

\smallskip

The action of the Galois group $\hat\Z^* \simeq \Gal(\Q^{ab}/\Q)$
on the coefficients of the $q$-expansion determines a homomorphism
\begin{equation}\label{cyclhom}
{\rm cycl}: \hat\Z^* \to \Aut(F).
\end{equation}

\smallskip

If $f$ is a continuous functions on $Z=\cR_2/\C^*$, we write
$$ f_{(g,\alpha)}(z) = f(g,\alpha,z) $$
so that $f_{(g,\alpha)}\in C(\H)$. For $p_N: M_2(\hat\Z)\to
M_2(\Z/N\Z)$ the canonical projection, we say that $f$ is of
level $N$ if
$$ f_{(g,\alpha)} = f_{(g,p_N(\alpha))} \ \ \ \ \forall (g,\alpha). $$
Then $f$ is completely determined by the functions
$$ f_{(g,m)} \in C(\H), \ \ \ \ \text{ for } m\in M_2(\Z/N\Z). $$

\smallskip

Notice that the invariance
$f(g\gamma,\alpha,z)=f(g,\gamma\alpha,\gamma(z))$, for all
$\gamma \in \Gamma$ and for all $(g,\alpha,z)\in \cU$, implies that
$f_{(g,m)| \gamma} = f_{(g,m)}$, for all $\gamma \in
\Gamma(N)\cap g^{-1}\Gamma g$, \ie $f$ is invariant under a congruence
subgroup.

\smallskip

The arithmetic algebra $\cA_{2,\Q}$ defined in \cite{CM} is a
subalgebra of continuous
functions on $Z=\cR_2/\C^*$ with the convolutions product
\eqref{Hecke-convol} and with the properties:
\begin{itemize}
\item The support of $f$ in $\Gamma\backslash\GL_2^+(\Q)$ is
finite.
\item The function $f$ is of finite level with
$$ f_{(g,m)} \in F \ \ \ \ \forall (g,m). $$
\item The function $f$ satisfies the {\em cyclotomic condition}:
$$ f_{(g,\alpha(u)m)} = {\rm cycl}(u) \, f_{(g,m)}, $$
for all $g\in \GL_2^+(\Q)$ diagonal and all $u\in \hat\Z^*$, with
$$ \alpha(u)=\begin{pmatrix} u& 0 \\ 0 & 1 \end{pmatrix} $$
and ${\rm cycl}$ as in \eqref{cyclhom}.
\end{itemize}

\smallskip

Here $F$ is the modular field, namely the field of modular
functions over $\Q^{ab}$. This is the union of the fields $F_N$ of
modular functions of level $N$ rational over the cyclotomic field
$\Q(\zeta_N)$.

\smallskip

The cycloomic condition is a consistency condition on the roots of
unity that appear in the coefficients of the $q$-series, which allows
for the existence of ``fabulous states'' (\cf \cite{CM}).

\smallskip

For $\alpha\in M_2(\hat\Z)$, let $G_\alpha\subset \GL_2^+(\Q)$ be
the set of
$$ G_\alpha = \{ g \in \GL_2^+(\Q): \,\, g\alpha \in M_2(\hat\Z)
\}. $$
Then, as shown in \cite{CM}, an element $y=(\alpha,z)\in
M_2(\hat\Z)\times \H$ determines a unitary representation of the Hecke
algebra $\cA$ on the Hilbert space $\ell^2(\Gamma\backslash
G_\alpha)$,
\begin{equation}\label{rep-y}
((\pi_y f)\xi)(g):= \sum_{s\in \Gamma\backslash G_\alpha}
f(gs^{-1},s\alpha, s(z))\, \xi(s), \ \ \ \ \forall g\in G_\alpha
\end{equation}
for $f\in \cA$ and $\xi\in \ell^2(\Gamma\backslash
G_\alpha)$.

\smallskip

Invertible $\Q$-lattices determine positive energy
representations, due to the fact
that the condition $g\alpha \in M_2(\hat \Z)$ for $g\in \GL_2^+(\Q)$ and
$\alpha\in \GL_2(\hat\Z)$ (invertible case) implies $g\in M_2(\Z)^+$,
hence the time evolution \eqref{GL2sigmat}
is implemented by the positive Hamiltonian with
spectrum $\{ \log \det (m) \}\subset [0,\infty)$ for $m\in
\Gamma\backslash M_2(\Z)^+$. The partition function of the
$\GL_2$-system is then $Z(\beta)=\zeta(\beta)\zeta(\beta-1)$, which
shows that one can expect the system to have two phase
transitions, which is in fact the case.

\smallskip

While the group $\GL_2(\hat\Z)$ acts by automorphisms of the
algebra of coordinates of $\cR_2/\C^*$, \ie the algebra of the
quotient \eqref{ShGL2NC}, the action of $\GL_2(\A_f)$ on the Hecke
algebra $\cA_2$ of coordinates of $\cR_2/\C^*$ is by {\em
endomorphisms}. More precisely, the group $\GL_2(\hat\Z)$ of deck
transformations of the modular tower still acts by automorphisms
on this algebra, while $\GL_2^+(\Q)$ acts by endomorphisms of the
$C^*$-dynamical system, with the diagonal $\Q^*$ acting by inner,
as in \eqref{innaction}.

\smallskip

The group of symmetries $\GL_2(\A_f)$ preserves the arithmetic
subalgebra. In fact, the group $\Q^*\backslash \GL_2(\A_f)$
has an important arithmetic meaning: a result of Shimura (\cf
\cite{Sh}, \cite{Lang}) characterizes the automorphisms of the modular
field by the exact sequence
\begin{equation}\label{shimuraexseq}
0 \to \Q^* \to \GL_2(\A_f) \to \Aut(F) \to 0.
\end{equation}

\smallskip

There is an induced action of $\Q^*\backslash \GL_2(\A_f)$
(symmetries modulo inner) on the KMS$_\beta$ states of
$(\cA,\sigma_t)$, for $\beta< \infty$. The action of
$\GL_2(\hat\Z)$ extends to KMS$_\infty$ states, while the action
of $\GL_2^+(\Q)$ on $\Sigma_\infty$ is defined by the action at
finite (large) $\beta$, by first ``warming up'' and then ``cooling
down''  as in \eqref{warming} (\cf \cite{CM}).

\smallskip

The result of \cite{CM} on the structure of KMS states for the $\GL_2$
system is as follows.

\smallskip

\begin{itemize}
\item There is no KMS state in the range $0<\beta\leq 1$.
\item In the range $\beta>2$ the set of extremal KMS states is
given by the invertible $\Q$-lattices, namely by the Shimura
variety $Sh(\GL_2,\H^\pm)$,
\begin{equation}\label{Ekmsbeta2}
\cE_\beta \cong \GL_2(\Q)\backslash \GL_2(\A) /\C^*.
\end{equation}
The explicit expression for these extremal KMS$_\beta$ states is
\begin{equation}\label{KMSexpr}
 \varphi_{\beta,L}(f)=\frac{1}{Z(\beta)} \sum_{m\in
\Gamma\backslash M_2^+(\Z)} f(1,m\alpha,m(z))\, \det(m)^{-\beta}
\end{equation}
where $L=(\alpha,z)$ is an invertible $\Q$-lattice.
\item At $\beta=\infty$, and for {\em generic} $L=(\alpha,\tau)$
invertible (where generic means $j(\tau)\notin\bar\Q$), the values
of the state $\varphi_{\infty,L}\in \cE_\infty$ on elements of
$\cA_{2,\Q}$ lie in an embedded image in $\C$ of the modular
field,
\begin{equation}\label{values}
\varphi(\cA_{2,\Q})\subset F_\tau,
\end{equation}
and there is an isomorphism
\begin{equation}\label{thetaphi1}
\theta_\varphi : \Aut_\Q(F_\tau)
\stackrel{\simeq}{\longrightarrow} \Q^* \backslash \GL_2(\A_f),
\end{equation}
depending on $L=(\alpha,\tau)$, which intertwines the Galois
action on the values of the state with the action of symmetries,
\begin{equation}\label{intertwineGL2}
\gamma\, \varphi(f) = \varphi( \theta_\varphi(\gamma) f), \ \ \ \
\forall f\in \cA_{2,\Q}, \ \ \forall\gamma\in \Aut_\Q(F_\tau).
\end{equation}
\end{itemize}

\section{Quantum statistical mechanics for imaginary quadratic
fields}\label{CMsect}

In the Kronecker--Weber case, the maximal abelian extension of
\(\Q\) is generated by the values  of the exponential
function at the torsion points \(\Q/{\Z}\) of the group \(\C/{\Z}
=\C^*\). Similarly, it is well known that the maximal abelian extension of
an imaginary quadratic field \(K\) is generated by the values of a certain
analytic function, the Weierstrass $\wp$-function, at the torsion points
$E_{tors}$ of an elliptic
curve and contains the \(j\)-invariant $j(E)$ of the elliptic curve.
Namely, in this case, the theory of complex multiplication of
elliptic curves provides a description of
 \(K^{ab}\). The ideal class group \({\rm Cl}(\mathcal O)\) is naturally
isomorphic to \({\rm Gal}(H/K)\), where \(H=K(j(E))\) is the
Hilbert class field of \(K\), \ie, its maximal abelian unramified
extension. In the case that \({\rm Cl}(\mathcal O)\) is trivial,
the situation of the CM case is exactly as for the field \(\Q\),
with \(\hat{\mathcal O}^*\) replacing \(\hat{\Z}^*\).

\smallskip

The construction of BC \cite{BC} was partially generalized to other
global fields (\ie number fields and function fields) in \cite{LvF,
hale,  Cohen}. The construction of \cite{hale} involves replacing the
ring $\O$ of integers of $K$ by a localized ring $\O_S$ which
is principal and then taking a cross product of the form
\begin{equation}\label{hl-alg}
C^*(K/\O_S) \rtimes \O_+^\times
\end{equation}
where $\O_+^\times$ is the sub semi-group
of $K^*$ generated by the generators of
prime ideals of $\O_S$. The symmetry group
is \(\hat{\mathcal O}_S^*\) and does not coincide
with what is needed
for class field theory except when the class number is $1$.
The construction of \cite{Cohen} involves a cross product
of the form
\begin{equation}\label{cohen-alg}
C^*(K/\O) \rtimes J^+
\end{equation}
where $J^+$ is a suitable adelic lift of the quotient
group $\J_K/\hat{\mathcal O}^*$. It gives the right
partition function namely the Dedekind zeta function
but not the expected symmetries.
The construction of \cite{LvF}
involves the algebra
\begin{equation}\label{KO-alg}
C^*(K/\O) \rtimes \O^\times \cong C(\hat\O) \rtimes
\O^\times
\end{equation}
and has symmetry group \(\hat{\mathcal O}^*\), while what is needed
for class field theory is a system with symmetry group \(\J_K/K^*\).
As one can see from the commutative diagram
\begin{eqnarray}
\diagram 1\rto & \hat\cO^*/\cO^* \rto\dto^{\simeq} & \J_K/K^*
\dto^{\simeq} \rto & {\rm Cl}(\cO)\dto^{\simeq} \rto & 1 \\
1\rto & \Gal(K^{ab}/H) \rto & \Gal(K^{ab}/K) \rto & \Gal(H/K) \rto & 1,
\enddiagram
\label{Hfielddiagr}
\end{eqnarray}
the action of \(\hat{\mathcal O}^*\) is sufficient only in
the case when the class number is one. In order to avoid the class number
one restriction in extending the results of \cite{BC} to imaginary
quadratic fields, it is natural to consider the universal situation:
the moduli space of elliptic curves with level structure, \ie, the
{\em modular tower}. Using the generalization of the BC case to
$\GL_2$ constructed in \cite{CM}, we will describe a system, which
does in fact have the right properties to recover the
explicit class field theory of imaginary quadratic fields from KMS
states. It is based on the geometric notions
of $K$-lattice and commensurability and extends to
quadratic fields the reinterpretation of the BC system
which was given in \cite{CM} in terms of $\Q$-lattices. The system we obtain 
in this section is, in fact, Morita equivalent to the one of \cite{LvF}. 
The two main new ingredients in our construction are the
choice of a natural rational subalgebra on which to evaluate the
KMS$_\infty$ states and the fact that the group of automorphisms
\(\hat{\mathcal O}^*/\O^*\) should be enriched by further symmetries, this
time given by {\em endomorphisms}, so that the actual group of
symmetries of the system is exactly \(\J_K/K^*\).

\medskip
\subsection*{$K$-lattices and commensurability}\hfill\medskip

In order to compare the BC system, the $\GL_2$ system and
the CM case, we give a definition of $K$-lattices, for $K$ an
imaginary quadratic field. The quantum statistical mechanical system
we shall construct to recover the explicit
class field theory of imaginary quadratic fields will be related to
commensurability of 1-dimensional $K$-lattices. This will be analogous
to the description of the BC system in terms of commensurability of
1-dimensional $\Q$-lattices. On the other hand, since there is a
forgetful map from 1-dimensional $K$-lattices to 2-dimensional
$\Q$-lattices, we will also be able to treat the CM case as a
specialization of the $\GL_2$ system at CM points.

\smallskip

Let \(\O = \Z + \Z \tau\) be the ring of integers of an imaginary
quadratic field $K=\Q(\tau)$; fix the imbedding \(K
\hookrightarrow \C\) so that \(\tau\in \H\).
Note that $\C$ then becomes a $K$-vector space
and in particular an $\O$-module.
The choice of $\tau$ as above also determines an
imbedding
\begin{equation}\label{qtau}
q_\tau: K \hookrightarrow M_2(\Q).
\end{equation}
The image of its restriction
$q_\tau: K^* \hookrightarrow \GL_2^+(\Q)$ is characterized by
the property that (\cf \cite{Sh}
Proposition 4.6)
\begin{equation}\label{fixtau}
q_\tau (K^*)=\{ g\in \GL_2^+(\Q): \,\,\, g(\tau)=\tau \}.
\end{equation}
For $g=q_\tau(x)$ with $x\in K^*$, we
have $\det(g)=\n(x)$, where $\n:K^*\to \Q^*$ is the norm map.

\smallskip

\begin{defn}\label{Klattices}
For $K$ an imaginary quadratic field, a 1-dimensional $K$-lattice
$(\Lambda,\phi)$ is a finitely generated $\O$-submodule $\Lambda\subset \C$,
such that $\Lambda\otimes_\O K \cong K$, together with a morphism of $\O$-modules
\begin{equation}\label{Kphi}
 \phi : K/\O \to K\Lambda/\Lambda.
\end{equation}
A 1-dimensional $K$-lattice is {\em invertible} if $\phi$ is an
isomorphism of $\O$-modules.
\end{defn}

\smallskip

Notice that in the definition we assume that the $\O$-module
structure is compatible with the embeddings of both $\O$ and
$\Lambda$ in $\C$.

\smallskip

\begin{lem}\label{1K2Q}
A 1-dimensional $K$-lattice is, in particular, a 2-dimensional $\Q$-lattice.
Moreover, as an $\O$-module, $\Lambda$ is projective.
\end{lem}

\proof First notice that $K\Lambda =\Q\Lambda$, since $\Q \O=K$. This, together with
$\Lambda\otimes_\O K \cong K$, shows that the $\Q$-vector space $\Q\Lambda$ is 2-dimensional.
Since $\R\Lambda =\C$, and $\Lambda$ is finitely generated as an abelian group,
this shows that $\Lambda$ is a lattice. The basis $\{ 1, \tau \}$ provides an identification
of $K/\O$ with $\Q^2/\Z^2$, so that we can view $\phi$ as a homomorphism of abelian groups
$\phi: \Q^2/\Z^2 \to \Q\Lambda/\Lambda$. The pair $(\Lambda,\phi)$
thus gives a two dimensional $\Q$-lattice.

As an $\O$-module $\Lambda$ is isomorphic to a finitely generated
$\O$-submodule of $K$, hence to an ideal in $\O$. Every ideal in a
Dedekind domain $\O$ is finitely generated projective over $\O$.

\endproof

\smallskip

The elliptic curve $E=\C/\Lambda$ has complex multiplication by $K$, namely
there is an isomorphism
\begin{equation}\label{iota}
\iota: K\stackrel{\simeq}{\to} \End(E)\otimes \Q.
\end{equation}
In general, for $\Lambda$ a lattice in $\C$, if the elliptic curve
$E=\C/\Lambda$ has complex multiplication (\ie there is an isomorphism
\eqref{iota} for $K$ an imaginary quadratic field), then the
endomorphisms of $E$ are given by $\End(E)=\O_\Lambda$, where
$\O_\Lambda$ is the order of $\Lambda$,
\begin{equation}\label{OLambda}
\O_\Lambda =\{ x\in K:\,\, x\Lambda \subset \Lambda \}.
\end{equation}
Notice that the elliptic curves $E=\C/\Lambda$, where $\Lambda$ is a
1-dimensional $K$-lattice, have $\O_\Lambda=\O$, the maximal order.
The number of distinct isomorphism classes
of elliptic curves $E$ with $\End(E)=\O$ is equal to the class number
$h_K$. All the other elliptic curves with complex multiplication by $K$
are obtained from these by isogenies.

\smallskip

\begin{defn}\label{Kcommens} Two 1-dimensional $K$-lattices
$(\Lambda_1,\phi_1)$ and
$(\Lambda_2,\phi_2)$ are commensurable if $K\Lambda_1=K\Lambda_2$
and $\phi_1=\phi_2$ modulo $\Lambda_1+\Lambda_2$.
\end{defn}

One checks as in the case of $\Q$-lattices (\cf \cite{CM}) that it is
an equivalence relation.

\smallskip

\begin{lem}\label{commens1K}
Two 1-dimensional $K$-lattices are commensurable iff
the underlying $\Q$-lattices are commensurable. Up to scaling, any
$K$-lattice $\Lambda$ is equivalent to a $K$-lattice $\Lambda'=\lambda
\Lambda \subset K \subset \C$. The latter is unique modulo $K^*$.
\end{lem}

\proof The first statement holds, since for 1-dimensional $K$-lattices
we have $K\Lambda =\Q\Lambda$.
For the second statement, the $K$-vector space
$K\Lambda$ is 1-dimensional. If $\xi$ is a generator, then
$\xi^{-1}\Lambda\subset K$. The remaining ambiguity is only by scaling
by elements in $K^*$.

\endproof

\smallskip

\begin{prop}\label{idclass}
For invertible 1-dimensional $K$-lattices, the element of $K_0(\O)$
associated to the $\O$-module $\Lambda$ is an invariant of the
commensurability class up to scaling.
\end{prop}

\proof Two invertible 1-dimensional $K$-lattices that are
commensurable are in fact equal. The same holds for lattices up to
scaling. Thus, the corresponding $\O$-module class is well defined.

\endproof

There is a canonical isomorphism $K_0(\O)\cong \Z + {\rm Cl}(\O)$ (\cf
Corollary 1.11, \cite{Milnor}), where the $\Z$ part is given by the
rank, which is equal to one in our case, hence the invariant of
Proposition \ref{idclass} is the class in ${\rm Cl}(\O)$.

\smallskip

In constrast to the Proposition above, every 1-dimensional $K$-lattice
is commensurable to a $K$-lattice whose $\O$-module structure is
trivial. This follows, since every ideal in $\O$ is commensurable to
$\O$.

\smallskip

\begin{prop}\label{1Klatrhos}
The data  $(\Lambda,\phi)$ of a 1-dimensional $K$-lattice are equivalent
to data $(\rho,s)$ of an element $\rho \in \hat\O$ and $s\in
\A^*_K/K^*$, modulo the action of $\hat\O^*$ given by $(\rho,s)\mapsto
(x^{-1}\rho, xs)$. Thus, the space of 1-dimensional $K$-lattices is
given by
\begin{equation}\label{space1dK}
\hat\O\times_{\hat\O^*} (\A^*_K/K^*).
\end{equation}
\end{prop}

\proof
The $\O$-module $\Lambda$ can be described in the form
$\Lambda_s=s_\infty^{-1}(s_f \hat\O \cap K)$, where
$s=(s_f,s_\infty)\in \A^*_K$. This satisfies $\Lambda_{ks}=\Lambda_s$ for all
$k\in (\hat\O^*\times 1)\,K^*\subset \A^*_K$. Indeed, up to scaling,
$\Lambda$ can be identified with 
an ideal in $\O$. These can be written in the form $s_f \hat\O \cap K$
(\cf \cite{Sh} \S 5.2). If $\Lambda_s=\Lambda_{s'}$, then $s_\infty'
s_\infty^{-1} \in K^*$ and one is reduced to the condition $s_f \hat\O
\cap K = s_f' \hat\O\cap K$, which implies $s_f' s_f^{-1} \in
\hat\O^*$. The data $\phi$ of the 1-dimensional $K$-lattice can be
described by the composite map $\phi = s_\infty^{-1}( s_f\circ \rho)$
\begin{eqnarray}
\diagram
K/\O \rto^\rho\dto^{\cong} & K/\O \rto\dto^{\cong} & K/(s_f \hat\O
\cap K)\dto^{\cong} \\
\A_{K,f}/\hat\O\rto^\rho & \A_{K,f}/\hat\O \rto^{s_f} & \A_{K,f}/s_f \hat\O
\enddiagram
\label{Kphi-rhos}
\end{eqnarray}
where $\rho$ is an element in $\hat\O$. By construction the map
$(\rho,s)\mapsto (\Lambda_s, s_\infty^{-1}( s_f\circ \rho))$
passes to the quotient
$\hat\O\times_{\hat\O^*} (\A^*_K/K^*)$ and the above shows
that it gives a bijection
with the space of 1-dimensional $K$-lattices.

\endproof

\smallskip

Notice that, even though $\Lambda$ and $\O$ are not isomorphic as $\O$-modules, 
on the quotients we have $K/\Lambda \simeq K/\O$ as $\O$-modules, with the 
isomorphism realized by $s_f$ in the diagram \eqref{Kphi-rhos}.

\smallskip

\begin{prop}\label{space1dKcomm}
Let $\A_K^\cdot =\A_{K,f}\times \C^*$ be the subset of ad\`eles of
$K$ with nontrivial archimedean component.
The map $\Theta(\rho,s)=\rho\, s$,
\begin{equation}\label{thetamap}
\Theta: \hat\O\times_{\hat\O^*} (\A^*_K/K^*) \to \A_K^\cdot/K^*,
\end{equation}
preserves commensurability and induces an identification
of the set of commensurability classes of 1-dimensional $K$-lattices
(not up to scale) with the space $\A_K^\cdot/K^*$.
\end{prop}

\proof The map is well defined, because $\rho s$ is invariant under
the action $(\rho,s)\mapsto (x^{-1}\rho, xs)$ of $x\in \hat\O^*$. It
is clearly surjective. It remains to show that two $K$-lattices
have the same image if and only if they are in the same
commensurability class. First we show that we can reduce to the case
of principal $K$-lattices, without changing the value of the map
$\Theta$. Given a $K$-lattice $(\Lambda,\phi)$, we write $\Lambda=
\lambda J$, where $J\subset \O$ is an ideal, hence
$\Lambda= \lambda( s_f\hat\O \cap K )$, where
$\lambda=s_\infty^{-1}\in \C^*$ and $s_f\in \hat\O\cap \A_{K,f}^*$.
Then $(\Lambda,\phi)$ is commensurate to the principal $K$-lattice
$(\lambda \O, \phi)$. If $(\rho,s)$ is the pair associated to
$(\Lambda,\phi)$, with $s=(s_f,s_\infty)$ as above, then
the corresponding pair $(\rho',s')$ for $(\lambda \O, \phi)$ is given by
$\rho'= s_f \, \rho$ and $s'=(1,s_\infty)$. Thus, we have
$\Theta(\Lambda,\phi)=\Theta(\lambda \O, \phi)$. We can then reduce to
proving the statement in the case of principal $K$-lattices
$(s_\infty^{-1} \O,s_\infty^{-1}\rho)$. In this case, the equality
$ s_\infty \rho = k s_\infty' \rho'$, for $k\in K^*$, means that we
have $s_\infty = k s_\infty'$ and $\rho =k \rho'$. In turn, this is
the relation of commensurability for principal $K$-lattices.

\endproof

\smallskip

Thus, we obtain, for 1-dimensional $K$-lattices, the following Lemma
as an immediate corollary,

\begin{lem} \label{Klatt-Shi}
The map defined as $\Upsilon: (\Lambda,\phi) \mapsto \rho\in
\hat\O/K^*$ for principal $K$-lattices extends to
an identification, given by $\Upsilon: (\Lambda,\phi) \mapsto
s_f\, \rho\in \A_{K,f}/K^*$, of the
set of commensurability classes of
1-dimensional $K$-lattices up to scaling  with the
quotient
\begin{equation}\label{AKDK2}
\hat\O/K^*=\A_{K,f}/K^*  .
\end{equation}
\end{lem}

\smallskip

The above quotient $\A_{K,f}/K^*$ has a description in terms of
elliptic curves, analogous to what we explained in the case of the
$\GL_2$-system. In fact, we can associate to $(\Lambda,\phi)$ the
data $(E,\eta)$ of an elliptic curve $E=\C/\Lambda$ with complex
multiplication \eqref{iota}, such that the embedding
$K\hookrightarrow \C$ determined by this identification and by the
action of $\End(E)$ on the tangent space of $E$ at the origin is
the embedding specified by $\tau$ (\cf \cite{Mi} p.28, \cite{Sh}
\S 5.1), and an $\A_{K,f}$-homomorphism
\begin{equation}\label{AKfhom}
\eta : \A_{K,f} \to \Lambda \otimes_\O \A_{K,f},
\end{equation}
The composite map
$$ \A_{K,f} \stackrel{\eta}{\to} \Lambda \otimes_\O \A_{K,f}
\stackrel{\simeq}{\to}  \A_{K,f}, $$
determines an element $s_f\,\rho\in \A_{K,f}$.
The set of equivalence classes of data $(E,\iota,\eta)$, where
equivalence is given by an isogeny of the elliptic curve compatible
with the other data, is the quotient $\A_{K,f}/K^*$.

\smallskip

This generalizes to the non-invertible case the analogous
result for invertible 1-dimensional $K$-lattices (data $(E,\iota,\eta)$,
with $\eta$ an isomorphism realized by an element $\rho \in
\A_{K,f}^*$) treated in \cite{Mi}, where the set of equivalence
classes is given by the id\`eles class group of the imaginary
quadratic field,
\begin{equation}\label{CKDK2}
\J_K/K^* = \GL_1(K)\backslash \GL_1(\A_{K,f})= C_K/D_K.
\end{equation}

\medskip
\subsection*{Algebras of coordinates}\hfill\medskip

We now describe the noncommutative algebra of coordinates of the
space of commensurability classes of 1-dimensional $K$-lattices up to
scaling.

\smallskip

To this purpose, we first consider
the
groupoid $\tilde\cR_K$
of the equivalence relation of commensurability on 1-dimensional
$K$-lattices (not up to scaling).
By construction, this groupoid is a subgroupoid of the groupoid
$\tilde \cR$ of commensurability classes of 2-dimensional
$\Q$-lattices. Its structure as a locally compact \'etale groupoid is
inherited from this embedding.

\smallskip

The groupoid $\tilde\cR_K$
 corresponds to the quotient $\A_K^\cdot/K^*$.
Its $C^*$-algebra is given up to Morita equivalence
by
the crossed product
\begin{equation}\label{crossAK}
C_0(\A_K^\cdot)\rtimes K^*.
\end{equation}

\smallskip

The case of commensurability classes of 1-dimensional
$K$-lattices up to scaling is more delicate. In fact, Lemma
\ref{Klatt-Shi} describes set theoretically the space of
commensurability classes of 1-dimensional
$K$-lattices up to scaling as the quotient $\A_{K,f}/K^*$. This has a
noncommutative algebra of coordinates, which is the crossed
product
\begin{equation}\label{crossAfK}
C_0(\A_{K,f}/\O^*)\rtimes K^*/\O^*.
\end{equation}
As we are going to show below, this is Morita equivalent to
the noncommutative algebra $\cA_K=C^*(G_K)$ obtained
by taking the quotient by scaling $G_K=\tilde\cR_K/\C^*$ of the groupoid
of the equivalence relation of commensurability.

\begin{prop} \label{qgroupoid} The quotient $G_K=
\tilde\cR_K/\C^*$ is a groupoid.
\end{prop}

\proof The simplest way to check this is to write $\tilde\cR_K$
as the union of the two groupoids $\tilde\cR_K= G_0 \cup G_1$
corresponding respectively to pairs of commensurable
$K$-lattices $(L,L')$ with $L=(\Lambda,0), L'=(\Lambda',0)$
and $(L,L')$ with $L=(\Lambda,\phi \neq 0), L'=(\Lambda',\phi' \neq 0)$.
The scaling action of $\C^*$ on $G_0$ is the identity on $\O^*$
and the corresponding action of $\C^*/\O^*$ is free on the units
of $G_0$. Thus the quotient $G_0/\C^*$ is a groupoid.
Similarly the action of $\C^*$ on $G_1$ is free
on the units
of $G_1$ and the quotient $G_1/\C^*$ is a groupoid.
\endproof

\smallskip

The quotient topology turns $G_K$ into a locally
compact \'etale groupoid.
 The algebra of coordinates
$\cA_K=C^*(G_K)$ is described equivalently
 by restricting the convolution product of
the algebra of $\tilde\cR_K$ to weight zero functions
with $\C^*$-compact support and then
passing to the $C^*$ completion, as in
the $\GL_2$-case.
In other words the algebra $\cA_K$ of the CM system is the convolution algebra of
weight zero functions on the groupoid $\tilde\cR_K$ of the equivalence
relation of commensurability on $K$-lattices.

\smallskip

Let us compare  the setting of \eqref{crossAfK}
\ie the groupoid $\A_{K,f}/\O^*\rtimes K^*/\O^*$
with the groupoid $G_K$.
In fact, the difference between these two settings can be seen by
looking at the case of $K$-lattices with $\phi=0$. In the first case,
this corresponds to the point $0\in \A_{K,f}/\O^*$, which has stabilizer
$K^*/\O^*$, hence we obtain the group $C^*$-algebra of $K^*/\O^*$. In the other
case, the corresponding groupoid is obtained as a quotient by $\C^*$
of the groupoid $\tilde\cR_{K,0}$ of pairs of commensurable
$\O$-modules (finitely generated of rank one) in $\C$.  In this case
the units of the groupoid $\tilde\cR_{K,0}/\C^*$ can be identified
with the elements of ${\rm Cl}(\O)$ and the reduced groupoid by any of
these units is the group $K^*/\O^*$. The general result below
gives the Morita equivalence in general.

\begin{prop}\label{RK0}
Let $\cH$ be the space of pairs
 up to scaling $((\Lambda,\phi),(\Lambda',\phi'))$  of commensurable
$K$-lattices, with $(\Lambda,\phi)$ principal. The space $\cH^{\,'}$ is
defined analogously with $(\Lambda',\phi')$ principal. The
corresponding bimodules give a Morita equivalence between the algebras
$C_0(\A_{K,f}/\O^*)\rtimes K^*/\O^*$ and $\cA_K=C^*(G_K)$.
\end{prop}

\proof The correspondence given by these bimodules has the effect of
reducing to the principal case. In that case the groupoid of the
equivalence relation (not up to scaling) is given by the crossed
product $\A_K^\cdot \rtimes K^*$. When taking the quotient by $\C^*$
we then obtain the groupoid $\A_{K,f}/\O^* \rtimes K^*/\O^*$.

\endproof

\smallskip

Recall that the $C^*$-algebra for the $\GL_2$-system is not unital,
the reason being that the space of 2-dimensional $\Q$-lattices up to
scaling is noncompact, due to the presence of the modulus $z\in \H$
of the lattice. When restricting to 1-dimensional $K$-lattices up to
scaling, this parameter $z$ affects only finitely many values,
corresponding to representatives $\Lambda=\Z + \Z z$ of the classes in
${\rm Cl}(\O)$. In fact one has,

\smallskip

\begin{lem}\label{unital}
The algebra $\cA_K$ is unital.
\end{lem}

\proof The set $G_K^{(0)}$  of units of
$G_K$ is the quotient of $\hat\O\times_{\hat\O^*}
(\A^*_K/K^*)$ by the action of $\C^*$. This gives the compact space
\begin{equation}\label{AKunits}
 X=G_K^{(0)}=\hat\O\times_{\hat\O^*} (\A^*_{K,f}/K^*).
\end{equation}
Notice that $\A^*_{K,f}/(K^*\times \hat\O^*)$ is ${\rm Cl}(\O)$. Since
the set of units is compact the convolution algebra is unital.

\endproof

\smallskip

By restriction from the $\GL_2$-system, there is a homomorphism $\n$ from
the groupoid $\tilde\cR_K$ to $\R_+^*$ given by the covolume of a
commensurable pair of $K$-lattices.
More precisely given such a pair $(L,L')=((\Lambda,\phi),(\Lambda',\phi'))$
we let
\begin{equation}\label{Kcov}
\vert L/L' \vert= {\rm covolume}(\Lambda')/{\rm covolume}(\Lambda)
\end{equation}

This is invariant under scaling
both lattices, so it is defined on $G_K=\tilde\cR_K/\C^*$. Up to scale, we
can identify the lattices in a commensurable pair with ideals in $\O$.
The covolume is then given by the ratio of the norms.
This defines a time evolution on the algebra $\cA_K$ by
\begin{equation}\label{timeK}
\sigma_t(f)(L,L')= \vert L/L' \vert^{it} f(L,L').
\end{equation}

\smallskip

We construct representations for the algebra $\cA_K$. For an \'etale groupoid
$G_K$, every unit $y\in G_K^{(0)}$ defines a
representation $\pi_y$ by left convolution of the algebra of
$G_K$ in the Hilbert space $\cH_y=\ell^2((G_K)_y)$,
where $(G_K)_y$ is the set of elements with source $y$.
The representations corresponding to points that have
a nontrivial automorphism group will no longer be irreducible.
As in the $\GL_2$-case, this defines the norm on $\cA_K$ as
\begin{equation}\label{normAK}
\| f \|= \sup_y \| \pi_y(f) \|.
\end{equation}

\smallskip

Notice that, unlike in the case of the $\GL_2$-system, we are dealing
here with amenable groupoids, hence the distinction between the
maximal and the reduced $C^*$-algebra does not arise.

\smallskip

\begin{lem}\label{posenLem}
\begin{enumerate}
\item
Given an invertible $K$-lattice $(\Lambda,\phi)$, the map
\begin{equation}\label{mapKinv}
(\Lambda',\phi')\mapsto J=\{ x\in \O | x\Lambda' \subset \Lambda \}
\end{equation}
gives a bijection of the set of $K$-lattices commensurable to
$(\Lambda,\phi)$ with the set of ideals in $\O$.
\item
Invertible $K$-lattices define positive energy
representations.
\item The partition function is the Dedekind zeta function.
\end{enumerate}
\end{lem}

\proof (1) As in Theorem 1.26 of \cite{CM} in the $\GL_2$-case, we use the fact
(Lemma 1.27 of \cite{CM}) that, if $\Lambda$ is an invertible 2-dimensional $\Q$-lattice and
$\Lambda'$ is commensurable to $\Lambda$, then $\Lambda\subset
\Lambda'$. The map above is well defined, since $J\subset \O$ is an
ideal. Moreover, $J\Lambda'=\Lambda$, since $\O$ is a Dedekind
domain. The map is injective, since $J$ determines $\Lambda'$ as the
$\O$-module $\{ x\in \C | x J \subset \Lambda \}$ and the
corresponding $\phi$ and $\phi'$ agree. This also shows that the map
is surjective. We then use the notation
\begin{equation}\label{Jinv}
J^{-1}(\Lambda,\phi) = (\Lambda',\phi).
\end{equation}

(2) For an invertible $K$-lattice, the above gives an identification
of $(G_K)_y$ with the set $\Jc$ of ideals $J\subset \O$. The
covolume is then given by the norm. The corresponding Hamiltonian is
of the form
\begin{equation}\label{KHamilt}
H\, \epsilon_J = \log \n(J)\,\, \epsilon_J.
\end{equation}

(3) The partition function of the CM system is then given by the
Dedekind zeta function
\begin{equation}\label{Zcm}
 Z(\beta)=\sum_{J \text{ ideal in } \O} \n(J)^{-\beta} =\zeta_K(\beta).
\end{equation}

\endproof

We give a more explicit description of the action
$L\mapsto J^{-1}L$ on $K$-lattices, which will be useful later. 

\begin{prop}\label{actJL}
Let $L=(\Lambda,\phi)$ be a $K$-lattices and $J\subset \O$ an
ideal. If $L$ is represented by a pair $(\rho,s)$, then $J^{-1}L$ is
represented by the commensurable pair $(s_J \rho,s_J^{-1}s)$, where
$s_J$ is a finite id\`ele such that $J=s_J \hat\O\cap K$.
\end{prop}

\proof The pair $(s_J \rho,s_J^{-1}s)$ defines an element in
$\hat\O\times_{\hat\O^*} \A_K^*/K^*$. In fact, first notice that $s_J$
is in fact in $\hat\O\cap \A_{K,f}^*$, hence the product $s_J \rho\in
\hat\O$. It is well defined modulo $\hat\O^*$, and by direct
inspection one sees that the class it defines in
$\hat\O\times_{\hat\O^*} \A_K^*/K^*$ is that of $J^{-1}L$.
By Proposition \ref{space1dKcomm} the $K$-lattice $(s_J
\rho,s_J^{-1}s)$ lies in the same commensurability class, since 
$\Theta(\rho,s)=\rho s=\Theta (s_J
\rho,s_J^{-1}s)$.

\endproof

\medskip
\subsection*{Symmetries}\hfill\medskip

We shall now adapt the discussion of symmetries of the $\GL(2)$ system
to $K$-lattices, adopting a contravariant  notation instead of the
covariant one used in \cite{CM}.

\begin{prop}\label{AKsymm}
The semigroup $\hat\O\cap \A_{K,f}^*$ acts on the algebra $\cA_K$ by
endomorphisms. The subgroup $\hat\O^*$ acts by
automorphisms. The subsemigroup $\O^\times$ acts by inner endomorphisms.
\end{prop}

\proof Given an ideal $J\subset \O$, consider the set of $K$-lattices
$(\Lambda,\phi)$ such that $\phi$ is well defined modulo
$J\Lambda$. We say, in this case, that the $K$-lattice
$(\Lambda,\phi)$ is divisible by $J$.
The above condition gives a closed and open subset of the set of
$K$-lattices up to scaling. We denote by $e_J\in \cA_K$ the
corresponding idempotent. Let $s\in \hat\O \cap \A_{K,f}^*$. Let
$J=s\hat\O \cap K$. Given a commensurable pair $(\Lambda,\phi)$ and
$(\Lambda',\phi')$, let
\begin{equation}\label{eJendo}
\theta_s( f) ((\Lambda,\phi),(\Lambda',\phi'))=\left\{
\begin{array}{ll} f((\Lambda,s^{-1}\phi),(\Lambda',s^{-1}\phi'))
  & \text{ both $K$-lattices are divisible by $J$} \\[2mm]
0 & \text{ otherwise. } \end{array} \right.
\end{equation}
The formula \eqref{eJendo} defines an endomorphism of
$\cA_K$ with range the algebra reduced by $e_J$. It is, in fact, an
isomorphism with the reduced algebra. Clearly, if $s\in \hat\O^*$ the
above defines an automorphism. If $s\in \O^\times$, the endomorphism
\eqref{eJendo} is inner. In fact, for $s\in \O^\times$, let $\mu_s\in
\cA_K$ be given by
\begin{equation}\label{mus}
 \mu_s((\Lambda,\phi),(\Lambda',\phi'))=\left\{
\begin{array}{ll} 1
  & \Lambda=s^{-1} \Lambda' \text{ and } \phi'=\phi \\[2mm]
0 & \text{ otherwise. } \end{array} \right.
\end{equation}
Then the range of $\mu_s$ is the projection $e_J$, where $J$ is the
principal ideal generated by $s$. Then we have
$$ \theta_s(f) =\mu_s \, f \, \mu_s^*,  \ \ \ \forall s\in \O^\times. $$

\endproof

The action of symmetries $\hat\O\cap \A_{K,f}^*$ is compatible with the time evolution,
$$ \theta_s \, \sigma_t = \sigma_t \, \theta_s, \ \ \ \ \ \forall
s\in\hat\O\cap \A_{K,f}^*, \, \forall t\in \R. $$
The isometries $\mu_s$ are eigenvectors of the time evolution, namely
$$ \sigma_t(\mu_s)=\n(s)^{it} \, \mu_s. $$
Thus, we can pass to the corresponding group of symmetries, modulo inner,
which is given by the id\`ele class group $\A_{K,f}^*/K^*$, which is
identified, by the class field theory isomorphism, with the Galois
group $\Gal(K^{ab}/K)$.

\smallskip

This shows that we have an action of the id\`ele class group
on the set of extremal KMS$_\beta$ states of the CM
system. The action of the subgroup $\hat\O^*/\O^*$ is by
automorphisms, while the action of the quotient group ${\rm Cl}(\O)$
is by endomorphisms, as we expected according to diagram \eqref{Hfielddiagr}.

In order to compute the value of KMS states
on the projection $e_J$ associated to an
ideal $J$ of the ring $\O$ of integers
(\ie the characteristic function of the
set of $K$-lattices divisible by $J$)
we introduce an isometry $\mu_J \in \cA_K$
such that its range is $e_J$. This isometry
is simply given with our notations by
\begin{equation}\label{muJ}
 \mu_J((\Lambda,\phi),(\Lambda',\phi'))=\left\{
\begin{array}{ll} 1
  & \Lambda=J^{-1} \Lambda' \text{ and } \phi'=\phi \\[2mm]
0 & \text{ otherwise. } \end{array} \right.
\end{equation}
which is similar to equation \eqref{mus} and reduces to that
one when $J$ is principal (generated by $s$).
Thus this would seem to imply that
it is not only the subsemigroup
$\O^\times$ that acts by inner but in fact a bigger one
using the isometries $\mu_J \in \cA_K$.
To see what happens one needs to compare the endomorphism
$f \to \mu_J f \mu_J^*$ with the endomorphism $\theta_s$.
In the first case one gets the following formula
\begin{equation}\label{eJendo1}
\mu_J f \mu_J^* ((\Lambda,\phi),(\Lambda',\phi'))=\left\{
\begin{array}{ll} f((J\,\Lambda,\phi),(J\,\Lambda',\phi'))
  & \text{ both $K$-lattices are divisible by $J$} \\[2mm]
0 & \text{ otherwise. } \end{array} \right.
\end{equation}
while in the second case it is given by formula \eqref{eJendo} \ie
\begin{equation}\label{eJendo2}
\theta_s( f) ((\Lambda,\phi),(\Lambda',\phi'))=\left\{
\begin{array}{ll} f((\Lambda,s^{-1}\phi),(\Lambda',s^{-1}\phi'))
  & \text{ both $K$-lattices are divisible by $J$} \\[2mm]
0 & \text{ otherwise. } \end{array} \right.
\end{equation}
Thus the key point is that the scaling is only allowed
by elements of $K^*$ and the scaling relation between the lattices
$(s\,\Lambda,\phi)$ and $(\Lambda,s^{-1}\phi)$ holds only
for $s\in K^*$ but not for ideles.
Thus even though the  $\mu_J$ always exist (for any ideal) they
implement the endomorphism $\theta_s$ only in the principal
case.

\medskip
\subsection*{Comparison with other systems}\hfill\medskip

\smallskip

It is also useful to see explicitly the relation of the algebra
$\cA_K$ of the CM system to the algebras
previously considered in generalizations of the Bost--Connes
results, especially those of \cite{Cohen}, \cite{hale}, and \cite{LvF}.
This will explain why the algebra $\cA_K$
contains exactly the amount of extra information to allow for
the full explict class field
theory to appear.

\smallskip

The partition function of the system considered in \cite{hale} 
agrees with the Dedekind zeta function only in the case of
class number one. A different system, which has partition function the
Dedekind zeta function in all cases, was introduced in
\cite{Cohen}. Our system also has as partition function the Dedekind
zeta function, independently of class number. It however differs from
the system of \cite{Cohen}. In fact, in the latter, which is a
semigroup crossed product, the natural
quotient of the $C^*$-algebra obtained by specializing at the fixed
point of the semigroup is the group ring of an extension of the class
group ${\rm Cl}(\O)$ by $K^*/\O^*$, while in our case, when
specializing similarly to the $K$-lattices with $\phi=0$, we obtain
an algebra Morita equivalent to the group ring of $K^*/\O^*$. Thus,
the two systems are not naturally Morita equivalent.

\smallskip

The system considered  in \cite{LvF} is analyzed there only under the
hypothesis of class number one. 
It can be recovered from our system, which has no restrictions on
class number, by reduction to those $K$-lattices that are
principal. Thus, the system of \cite{LvF} is Morita equivalent to
our system (\cf Proposition \ref{RK0}). 

\smallskip

Notice, moreover, that the crossed product algebra $C(\hat\O)\rtimes
\O^\times$ considered in some generalizations of the BC system is more
similar to the ``determinant part'' of the
$\GL_2$-system (\cf Section 1.7 of \cite{CM}), namely to the algebra
$C(M_2(\hat\Z)) \rtimes M_2^+(\Z)$, than to the full $\GL_2$-system.

\section{KMS states and class field theory for imaginary quadratic
fields}\label{KMSsect}

The relation between the CM and the $\GL_2$-system provides us with a
choice of an arithmetic subalgebra $\cA_{K,\Q}$ of $\cA_K$. This is
obtained by restricting elements of $\cA_{2,\Q}$ to the
$\C^*$-quotient $G_K$ of the subgroupoid $\tilde\cR_K\subset \tilde\cR_2$.
Notice that, for the CM system, $\cA_{K,\Q}$ is a subalgebra of
$\cA_K$, not just a subalgebra of unbounded multipliers as in the
$\GL_2$-system, because of the fact that $\cA_K$ is unital.

\smallskip

We are now ready to state the main result on the CM case. The
following theorem will be proved in various steps in this section.

\smallskip

\begin{thm}\label{KMSCM}
Consider the system $(\cA_K,\sigma_t)$ described in the previous
section. The extremal KMS states of this system satisfy:
\begin{itemize}
\item In the range $0<\beta\leq 1$ there is a unique KMS state.
\item For $\beta>1$, extremal KMS$_\beta$ states are parameterized by
invertible $K$-lattices,
\begin{equation}\label{KMSbeta1}
\cE_\beta \simeq \A_{K,f}^*/K^*
\end{equation}
with a free and transitive action of the id\`ele class group of $K$ as
symmetries.
\item In this range, the extremal KMS$_\beta$ state associated to an
invertible $K$-lattice $L=(\Lambda,\phi)$ is of the form
\begin{equation}\label{varphibetaL}
\varphi_{\beta,L}(f) = \zeta_K(\beta)^{-1} \,\, \sum_{J \text{ ideal
in }
\O} f(J^{-1}L ,J^{-1} L)\,\, \n(J)^{-\beta},
\end{equation}
where $\zeta_K(\beta)$ is the Dedekind zeta function, and $J^{-1}L$
defined as in \eqref{Jinv}.
\item The set of extremal KMS$_\infty$ states (as weak limits of
KMS$_\beta$ states) is still given by \eqref{KMSbeta1}.
\item The extremal KMS$_\infty$ states $\varphi_{\infty,L}$ of the CM
system, evaluated on the arithmetic subalgebra $\cA_{K,\Q}$, take
values in $K^{ab}$. The class field theory isomorphism \eqref{CFTiso}
intertwines
the action of $\A_{K,f}^*/K^*$ by symmetries of the system
$(\cA_K,\sigma_t)$ and the action of $\Gal(K^{ab}/K)$ on the image
of $\cA_{K,\Q}$ under the extremal KMS$_\infty$ states,
\begin{equation}\label{intertwGalK}
 \alpha \circ \varphi_{\infty,L} = \varphi_{\infty,L} \circ
\theta^{-1}(\alpha), \ \ \ \forall \alpha \in \Gal(K^{ab}/K).
\end{equation}
\end{itemize}
\end{thm}

\smallskip

The proof of this statement is given in the following subsections.

\smallskip

Notice that the result stated above is substantially different from
the $\GL_2$-system. This is not surprising, as the following general
fact illustrates.
Given an \'etale  groupoid $\cG$ and a full subgroupoid
$\cG'\subset \cG$, let $\rho$ be a homomorphism $\rho: \cG \to
\R^*_+$. The inclusion $\cG'\subset \cG$ gives a correspondence
between the $C^*$-algebras associated to $\cG'$ and $\cG$, compatible
with the time evolution associated to $\rho$ and its restriction to
$\cG'$. The following simple example, however, shows that, in general,
the KMS states for the $\cG'$ system do not map to KMS states for the
$\cG$ system. We let $\cG$ be the groupoid with units $\cG^{(0)}$
given by an infinite countable set, and morphisms given by all pairs
of units. Consider a finite subset of $\cG^{(0)}$ and let $\cG'$ be
the reduced groupoid. Finally, let $\rho$ be trivial. Clearly, the
$\cG'$ system admits a KMS state for all temperatures given by the
trace, while, since there is no tracial state on the compact
operators, the $\cG$ system has no KMS states.

\medskip
\subsection*{KMS states at low temperature}\hfill\medskip

The partition function $Z_K(\beta)$ of \eqref{Zcm} converges for
$\beta>1$. We have also seen in the previous section that invertible
$K$-lattices $L=(\Lambda,\phi)$ determine positive energy
representations of $\cA_K$ on the Hilbert space
$\cH=\ell^2(\Jc)$ where $\Jc$ is the set of
ideals of $\Oc$. Thus, the formula
\begin{equation}\label{CM-KMS-Tr}
\varphi_{\beta,L}(f)=\frac{\Tr\left( \pi_L(f)\, \exp(-\beta H)
\right)}{\Tr(\exp(-\beta H))}
\end{equation}
defines an extremal KMS$_\beta$ state, with the Hamiltonian $H$ of
\eqref{KHamilt}. These states are of the form \eqref{varphibetaL}.
It is not hard to see that distinct elements in $\A_{K,f}/K^*$ define
distinct states $\varphi_{\beta,L}$.

\smallskip

This shows that we have an injection of $\A_{K,f}/K^*\subset
\cE_\beta$. We need to show that, conversely, every extremal KMS$_\beta$
state is of the form \eqref{varphibetaL}.

\smallskip

In order to prove the second and third statements
 of Theorem \ref{KMSCM} we shall proceed in
two steps. The first
shows (Proposition \ref{measure} below) that KMS$_\beta$ states
are given by measures on the space $X$ of $K$-lattices (up to
scaling). The second shows that when $\beta > 1$ this measure is
carried by the commensurability classes of invertible
$K$-lattices.

\smallskip

 We first describe the elements $\g \in G_K$
such that $s(\g)=r(\g)$ \ie $\g=(L,L')\in\tilde \cR_K$
such that the classes of $L$ and $L'$ modulo scaling are the
same elements of $X=G_K^{(0)}$.
Modulo scaling we can assume that the lattice
$\Lambda \subset K$ and since $L$ and $L'$ are commensurable
it follows that $\Lambda' \subset K$. Then by hypothesis there
exists $\lambda \in \C^*$ such that $\lambda L=L'$.
One has $\lambda \in K^*$ and $\phi'= \lambda \phi$
but by commensurability of the pair one also has
$\phi'= \phi$ modulo $ \Lambda+\Lambda'$.
Writing $\lambda=\frac{a}{b}$ with $a,b \in \O$ we get
$a \phi= b\phi$ and $\lambda=1$ unless $\phi =0$.
We thus get, with
$p$ the projection from $K$-lattices to their class
$p(L)\in X$ modulo scaling,

\smallskip

\begin{lem} \label{fixedpt} Let $\g=(L,L')\in\tilde \cR_K$
with $p(L)=p(L')\in X$. Then either $L=L'$
or $\phi=\phi'=0$.
\end{lem}

 We let $F$ be the finite closed
subset of $X$ given by the set of $K$-lattices
up to scaling such that $\phi=0$. Its cardinality is the
class number of $K$. The groupoid $G_K$ is the union
$G_K =G_0 \cup G_1$ of the reduced groupoids by $F\subset X$
and its complement.

\begin{lem} \label{corresp} Let   $\g\in G_1
\backslash G_1^{(0)}$.
There exists  a neighborhood $V$ of $\g$ in $G_K$  such that
$$
 r(V)\cap  s(V)=\, \emptyset
$$
where $r$ and $s$ are the range and source maps of  $G_K$.
\end{lem}

\proof Let $\g$ be the class modulo
scaling of the commensurable pair
$(L,L')$. By Lemma \ref{fixedpt}
one has $p(L)\neq p(L')\in X$. If the classes
of $\Lambda$ and $\Lambda'$ in $K_0(\O)$ are different
one just takes $V$ so that all elements
$\g_1=(L_1,L'_1)\in V$ are in the corresponding
classes which ensures $ r(V)\cap  s(V)=\, \emptyset$.
Otherwise there exists $\lambda \in K^*$
such that  $\Lambda'=\lambda \Lambda$
and since  $\Lambda'\neq \Lambda$ one has $\lambda\notin \O^*$.
One has $\phi'=\phi \neq 0$. Thus one is reduced to showing that
given $\rho \in \hat\O$, $\rho \neq 0$, and $\lambda \in K^*$,
$\lambda\notin \O^*$ there exists a neighborhood $W$
of $\rho $ in $\hat\O$ such that $\lambda W \cap \O^*W= \emptyset
$. This follows using a place $v$ such that $\rho_v\neq 0$, one
has $\lambda \rho_v \notin \O^*\rho_v$ and the same holds in a
suitable neighborhood.
\endproof

\smallskip

 We can now prove the following.
\begin{prop} \label{measure}
 Let $\beta >0$ and $\varphi$ a  KMS$_\beta$
state  on $(\cA_K,\sigma_t)$. Then there exists a probability measure
$\mu$ on $X$ such that
$$
\varphi(f)=\,\int_X\,f(L,L)\,d\mu(L)\qqq f\in \cA_K\,.
$$
\end{prop}

\proof It is enough to show
that $\varphi(f)=0$ provided $f$ is a continuous
 function with compact support
 on $G_K$ with support disjoint from
$G_K^{(0)}$.
 Let $h_n\in C(X),\,0\leq h_n\leq 1$ with
support disjoint from $F$ and converging pointwise to $1$ in the
complement of $F$. Let $u_n \in \cA_K$ be supported by the
diagonal and given by $h_n$ there.
\smallskip

 The formula
\begin{equation}\label{Phif}
\Phi(f)(\Lambda,\Lambda'):= f((\Lambda,0),(\Lambda',0))) \ \ \ \forall f\in \cA_K
\end{equation}
defines a homomorphism of $(\cA_K,\sigma_t)$ to the $C^\ast$ dynamical
system $(C^*(G_0),\sigma_t)$ obtained by specialization to
pairs of $K$-lattices with $\phi=0$ as in \cite{CM}.

\smallskip

 Since there are unitary eigenvectors for $\sigma_t$ for
non trivial eigenvalues in the system
$(C^*(G_0),\sigma_t)$ it has no non-zero KMS$_\beta$ positive functional.
This shows that the pushforward
of $\varphi$ by $\Phi$ vanishes and
by Proposition 5 of \cite{CM}  that, with the
notation introduced above,
$$
\varphi(f)= \lim_n \,\varphi(f * u_n)\,.
$$
Thus, since $(f * u_n)(\g)=\,f(\g)\,h_n(s(\g))$, we can assume that
$f(\g)=0$ unless $s(\g)\in C$, where $C\subset X$ is a compact
subset disjoint from $F$. Let $L\in C$ and $V$ as in lemma
\ref{corresp} and let $h\in C_c(V)$. Then, upon applying the
KMS$_\beta$ condition  to the pair  $a,b$ with $a=f$
and $b$ supported by the diagonal and equal to $h$ there.
One gets $\varphi(b*f)=\varphi(f*b)$. One has
$(b*f)(\g)=\,h(r(\g))\,f(\g)$. Applying this to $f*b$ instead of
$f$ and using $h(r(\g))\,h(s(\g))=0 \qqq \g\in V$, we get
$\varphi(f*b^2)=0$ and $\varphi(f)=0$, using a partition of unity.
\endproof

\smallskip

\begin{lem} \label{proba}
Let $\varphi$ be a  KMS$_\beta$
state  on $(\cA_K,\sigma_t)$. Then, for any ideal $J\subset \O$
one has
$$
\varphi(e_J)=\,\n(J)^{-\beta}\,.
$$
\end{lem}

\proof For each ideal $J$ we let $\mu_J\in
\cA_K$ be given as above by \eqref{muJ}
$$
  \mu_J((\Lambda,\phi),(\Lambda',\phi'))=\left\{
\begin{array}{ll} 1
  & \Lambda=J^{-1} \Lambda' \text{ and } \phi'=\phi \\[2mm]
0 & \text{ otherwise. } \end{array} \right.
$$
One has $\sigma_t(\mu_J)= \n(J)^{it}\,\mu_J \forall t\in \R$
while $\mu_J^*\,\mu_J=1$ and $\mu_J\,\mu_J^*=\,e_J$ thus the
answer follows from the KMS condition. 
\endproof

\smallskip

Given Proposition \ref{1Klatrhos} above, we make the following
definition.

\begin{defn}\label{qinvKlat}
A $K$-lattice is {\em quasi-invertible} if the $\rho$ in
Proposition \ref{1Klatrhos} is in $\hat\O \cap \A^*_{K,f}$.
\end{defn}

Then we have the following result.

\begin{lem}\label{invidele}
\begin{enumerate}
\item A $K$-lattice $(\Lambda,\phi)$ that is divisible by only
finitely many ideals is either quasi-invertible, or there is a
finite place $v$ such that $\phi_v=0$.
\item A quasi-invertible $K$-lattice is
commensurable to a unique invertible $K$-lattice.
\end{enumerate}
\end{lem}

\proof Let $(\rho,s)$ be associated to the $K$-lattice
$(\Lambda,\phi)$ as in Proposition \ref{1Klatrhos}. If $\rho\notin
\A^*_{K,f}$, then either there exists a place $v$ such that
$\rho_v=0$, or $\rho_v\neq 0$ for all $v$ and there exists
infinitely many places $w$ such that $\rho_w^{-1}\notin \O_w$,
where $\O_w$ is the local ring at $w$. This shows that the
$K$-lattice is divisible by infinitely many ideals. For the second
statement, if we have $\rho\in \A^*_{K,f}$, we can write it as a
product $\rho=s_f' \rho'$ where $\rho'=1$ and $s_f'=\rho$. The
$K$-lattice obtained this way is commensurable to the given one by
Proposition \ref{space1dKcomm} and is invertible.

\endproof

\smallskip

Let us now complete the proof of the second and third
statements  of Theorem \ref{KMSCM}.
Let $\varphi$ be a KMS$_\beta$ state.
Proposition \ref{measure} shows that there is a
probability measure $\mu$ on $X$
such that
$$
\varphi(f)=\,\int_X\,f(L,L)\,d\mu(L)\qqq f\in \cA_K\,.
$$
With $L=(\Lambda, \phi) \in X$, Lemma \ref{proba} shows that the
probability $\varphi(e_J)$
that an ideal $J$ divides $L$ is $\n(J)^{-\beta}$. Since
 the series $\sum
\n(J)^{-\beta}$ converges for $\beta>1$, it follows (\cf
\cite{Rudin} Thm. 1.41) that, for almost all $L\in X$, $L$ is only
divisible by a finite number of ideals. Notice that the KMS
condition implies that the measure defined above gives measure
zero to the set of $K$-lattices $(\Lambda,\phi)$ such that
$\phi_v=0$ for some finite place $v$.

By the first part of Lemma \ref{invidele}, the measure $\mu$ gives
measure one to quasi-invertible $K$-lattices. (Notice that these
$K$-lattices form a Borel subset which is not closed.) Then, by
the second part of Lemma \ref{invidele}, the KMS$_\beta$ condition
shows that the measure $\mu$ is entirely determined by its
restriction to invertible $K$-lattices, so that, for some
probability measure $\nu$,
$$
\varphi=\, \int \,\varphi_{\beta,L}\,\,d\nu(L).
$$
It follows that the Choquet simplex of extremal KMS$_\beta$ states
is the space of probability measures on the
compact space
of invertible $K$-lattices modulo scaling
and its extreme points are the $\varphi_{\beta,L}$. $\Box$\\

The action of the symmetry group $\J_K/K^*$ on $\cE_\beta$ is then
free and transitive. In fact, recall that (Lemma 1.28 \cite{CM}) the
action of $\GL_2^+(\Q)$ on 2-dimensional $\Q$-lattices has as
its only fixed points the $\Q$-lattices $(\Lambda,\phi)$
with $\phi=0$.

\smallskip

The action is given explicitly, for $L=(\Lambda,\phi)$
an invertible $K$-lattice
and for $s\in \hat\O^*/\O^*\subset \J_K/K^*$, by
\begin{equation}\label{actionstates1}
(\varphi_{\beta,L}\circ \theta_s) (f)= Z_K(\beta)^{-1} \,\, \sum_{J\in
\Jc} f((J^{-1}\Lambda,s^{-1}\phi),(J^{-1}\Lambda,s^{-1}\phi)) \,\, \n(J)^{-\beta}
=\varphi_{\beta,(\Lambda,s^{-1}\phi)}(f),
\end{equation}
More generally, for $s\in \hat\O \cap \A_{K,f}^*$ let
$J_s=s\hat\O \cap K$, one then has,
\begin{equation}\label{actionstates2}
\begin{array}{rl}
(\varphi_{\beta,L}\circ \theta_s) (f)= & Z_K(\beta)^{-1} \sum_{J\in
\Jc} \theta_s(f)((J^{-1}\Lambda,\phi),(J^{-1}\Lambda,\phi))
 \,\, \n(J)^{-\beta} \\[3mm]
= & Z_K(\beta)^{-1} \sum_{J\supset J_s}
f((J^{-1}\Lambda,s^{-1}\phi),(J^{-1}\Lambda,s^{-1}\phi)) \,\,
\n(J)^{-\beta}\\[3mm]
= & Z_K(\beta)^{-1} \sum_{J\in \Jc}
f(J^{-1}L_s,J^{-1}L_s) \,\,
\n(J \,J_s)^{-\beta} = \n(J_s)^{-\beta}\,\varphi_{\beta,L_s}(f), \end{array}
\end{equation}
where $L_s$ is the invertible
$K$-lattice $(J_s^{-1}\Lambda,s^{-1}\phi)$.
To prove the last equality one uses the basic
property of Dedekind rings that any ideal $J\supset J_s$
can be written as a product $J=\,J'\,J_s$.

\medskip
\subsection*{KMS states at zero temperature and Galois
action}\hfill\medskip

The weak limits as $\beta\to \infty$ of states in $\cE_\beta$ define
states in $\cE_\infty$ of the form
\begin{equation}\label{Kground}
 \varphi_{\infty,L}(f) = f(L,L).
\end{equation}

\smallskip

Some care is needed in defining the action of the symmetry group
$\A_{K,f}/K^*$ on extremal states at zero temperature. In fact, as it
happens also in the $\GL_2$-case, for an invertible $K$-lattice
 evaluating $\varphi_{\infty,L}$ on
$\theta_s(f)$ does not give a nontrivial action in the case of
endomorphisms. However, there is a nontrivial action induced on
$\cE_\infty$ by the action on $\cE_\beta$ for finite $\beta$ and it is
obtained as
\begin{equation}\label{warmcool}
 \Theta_s(\varphi_{\infty,L})(f) = \lim_{\beta\to \infty} \left(
W_\beta(\varphi_{\infty,L})\circ \theta_s \right)(f),
\end{equation}
where $W_\beta$ is the ``warm up'' map \eqref{warming}. This gives
\begin{equation}\label{Thetas}
\Theta_s(\varphi_{\infty,L}) = \varphi_{\infty,L_s}.
\end{equation}

\smallskip

Thus the action of the symmetry group $\J_K/K^*$ is given by
\begin{equation}\label{actsym}
L=\,(\Lambda,\phi) \to L_s=\,(J_s^{-1}
\Lambda,s^{-1}\phi)\qqq s\in \J_K/K^*\,.
\end{equation}

When we evaluate states $\varphi_{\infty,L}$ on elements $f\in
\cA_{K,\Q}$ of the arithmetic subalgebra we obtain
\begin{equation}\label{valuesKab}
\varphi_{\infty,L}(f)=f(L,L)=g(L)\,,
\end{equation}
where the function $g$ is the lattice function of weight $0$
obtained as the restriction of $f$ to the diagonal. By
construction of $\cA_{K,\Q}$, one obtains in this way all the
evaluations $f\mapsto f(z)$ of elements of the modular field $F$
on the finitely many modules $z \in \H$ of the classes of
$K$-lattices.

\smallskip

The modular functions \(f\in F\) that are defined at \( \tau \)
define a subring \(B\) of \(F\). The theory of complex
multiplication (\cf \cite{Sh}) shows that the subfield
$F_\tau\subset \C$ generated by the values $f(\tau)$, for $f\in
B$, is the maximal abelian extension of $K$ (we have fixed an
embedding $K\subset \C$),
\begin{equation}\label{FtauKab}
F_\tau = K^{ab}.
\end{equation}

\smallskip

Moreover, the action of $\alpha\in\Gal(K^{ab}/K)$ on the values
$f(z)$ is given by
\begin{equation}\label{Shirec}
\alpha f(z) = f^{\,\sigma q_\tau \theta^{-1}(\alpha)}(z).
\end{equation}
In this formula the notation $f\mapsto f^\gamma$ denotes the action of
an element $\gamma\in \Aut(F)$ on the elements $f\in F$, the map
$\theta$ is the class field theory isomorphism
\eqref{CFTiso}, $q_\tau$ is the embedding of $\A_{K,f}^*$ in
$\GL_2(\A)$ and $\sigma$ is as in the diagram with exact rows
\begin{eqnarray}
\diagram 1\rto & K^* \rto^{\iota\qquad} & \GL_1(\A_{K,f})
\dto^{q_\tau} \rto^{\theta} &
\Gal(K^{ab}/K) \rto & 1 \\
1\rto & \Q^* \rto & \GL_2(\A_f) \rto^{\sigma} & \Aut(F) \rto & 1.
\enddiagram
\label{Shirecdiagr}
\end{eqnarray}

\smallskip

Thus, when we act by an element $\alpha\in\Gal(K^{ab}/K)$ on the
values on $\cA_{K,\Q}$ of an extremal KMS$_\infty$ state we have
\begin{equation}\label{Galphi}
\alpha\, \varphi_{\infty,L}(f) = \varphi_{\infty,L_s}(f)
\end{equation}
where $s=\theta^{-1}(\alpha)\in \J_K/K^*$.

\smallskip

This corresponds to the result of Theorem 1.39 of \cite{CM} for the
case of 2-dimensional $\Q$-lattices (see
equations (1.130) and following in \cite{CM})
with the slight nuance that we used there a
covariant notation for the Galois action rather than
the traditional contravariant one $f\mapsto f^\gamma$.

\medskip
\subsection*{Uniqueness of high temperature KMS state}\hfill\medskip

The proof follows along the line of \cite{nes}. We first discuss
uniqueness. By Proposition \ref{measure}, one obtains a measure
$\mu$ on the space $X$ of $K$-lattices up to scale. As in Lemma
\ref{proba}, this measure fulfills the quasi-invariance condition
\begin{equation}\label{qinvcond}
\int_X \mu_J \,f \, \mu_J^* \,\, d\mu = \n(J)^{-\beta}  \, \int_X
f \, d\mu,
\end{equation}
for all ideals $J$, where $\mu_J$ is as in \eqref{muJ}. To prove
uniqueness of such a measure, for $\beta\in (0,1]$, one proceeds
in the same way as in \cite{nes}, reducing the whole argument to
an explicit formula for the orthogonal projection $P$ from
$L^2(X,d\mu)$ to the subspace of functions invariant under the
semigroup action
\begin{equation}\label{JL}
 L=(\Lambda,\phi) \mapsto J^{-1}L,
\end{equation}
which preserves commensurability. As in \cite{nes}, one can obtain
such formula as a weak limit of the orthogonal projections $P_A$
associated to finite sets $A$ of non-archimedean places.

\smallskip

Let $A$ be a finite set of non-archimedean places. Let $\sJ_A$ be
the subsemigroup of the semigroup $\sJ$ of ideals, generated by
the prime ideals in $A$. 
 Any element $J\in \sJ_A$ can be uniquely written as a
product
\begin{equation}\label{prodJ}
J= \prod_{v\in A} J_v^{n_v},
\end{equation}
where $J_v$ is the prime ideal associated to the place $v\in A$.

\smallskip

\begin{lem}\label{partitionA}
Let $L=(\Lambda,\phi)$ be a $K$-lattice such that $\phi_v\neq 0$
for all $v\in A$. Let $J\in \sJ_A$, 
$J= \prod_{v\in A} J_v^{n_v}$ be the smallest ideal dividing
$L$. Let $(\rho,s)\in \hat\O\times_{\hat\O^*} \A^*_K/K^*$ be the
pair associated to $L$. Then, for each $v\in A$, the valuation of
$\rho_v$ is equal to $n_v$.
\end{lem}

\proof Let $(\rho,s)$ be as above, and $m_v$ be the valuation of
$\rho_v$. Then it is enough to show that an ideal $J$ divides $L$
if and only if $J$ is of the form \eqref{prodJ}, with $n_v\leq
m_v$. The map $\phi$ is the composite of multiplication by $\rho$
and an isomorphism, as in the diagram \eqref{Kphi-rhos}, hence the
divisibility is determined by the valuations of $\rho_v$.
\endproof

\smallskip

\begin{defn}\label{Ainv}
With $A$ as above we shall say that a $K$-lattice 
$L=(\Lambda,\phi)$ is $A$-invertible iff the valuation of 
$\rho_v$ is equal to zero far all $v \in A$.
\end{defn}

\smallskip

We now define basic test functions associated to a Hecke
Grossencharakter. Given such a character $\chi$, the restriction
of $\chi$ to $\hat\O^*$ only depends on the projection on
$\hat\O^*_{B_\chi} =\prod_{v\in B_\chi} \hat\O^*_v$, for $B_\chi$
a finite set of non-archimedean places. Let $B$ be a finite set of
non-archimedean places $B\supset B_\chi$. We consider the function
$f=f_{B,\chi}$ on $\hat\O\times_{\hat\O^*} \A^*_K/K^*$, which is
obtained as follows. For $(\rho,s)\in \hat\O\times_{\hat\O^*}
\A^*_K/K^*$, we let $f=0$ unless $\rho_v\in \hat\O^*_v$ for all
$v\in B$, while $f(\rho,s)=\chi(\rho' s)$, for any $\rho'\in
\hat\O^*$ such that $\rho'_v =\rho_v$ for all $v\in B$. This is
well defined, since the ambiguity in the choice of $\rho'$ does
not affect the value of $\chi$, since $B_\chi \subset B$. The
function obtained this way is continuous.

\smallskip

Let $H(B)$ be the subspace of $L^2(X,d\mu)$ of functions that only
depend on $s$ and on the projection of $\rho$ on $\hat\O_B$. Let
us consider the map that assigns to a $K$-lattice $L$ the smallest
ideal $J\in \sJ_B$ dividing $L$, extended by zero if some
$\phi_v=0$. By Lemma \ref{partitionA}, the value of this map only
depends on the projection of $\rho$ on $\hat\O_B$. By construction
the corresponding projections $E_{B,J}$ give a partition of unity
on the Hilbert space $H(B)$. Note that $E_{B,0}=0$, since the
measure $\mu$ gives measure zero to the set of $K$-lattices with
$\phi_v=0$ for some $v$.

\smallskip

Let $V_J f (L)=f(J^{-1}L)$ implementing the semigroup action
\eqref{JL}. For $J\in \sJ_B$, the operator $\n(J)^{-\beta/2}
V_J^*$ maps isometrically the range of $E_{B,\O}$ to the range of
$E_{B,J}$.

\smallskip

\begin{lem}\label{VJchi}
The functions $V_J^* f_{B,\chi}$ span a dense subspace of $H(B)$.
\end{lem}

\proof It is sufficient to prove that the $f_{B,\chi}$ form a
dense subspace of the range of $E_{B,\O}$. The image of $\hat\O^*$
in $\hat\O^*_B \times \A^*_K/K^*$ by the map $u \mapsto
(u,u^{-1})$ is a closed normal subgroup. We let $\hat\O^*_B
\times_{\hat\O^*} \A^*_K/K^*$ be the quotient. This is a locally
compact group. The quotient $G_B$ by the connected component of
identity $D_K$ in $\A^*_K/K^*$ is a compact group. Then $C(G_B)$
is identified with a dense subspace of the range of $E_{B,\O}$.
The characters of $G_B$ are the Grossencharakters $\chi$ that
vanish on the connected component of identity and such that
$B_\chi\subset B$. Thus, by Fourier transform, we obtain the
density result.
\endproof

\smallskip

Let $A$ be a finite set of non-archimedean places, and $\sJ_A$ as
above.
Let $H_A$ be the subspace of functions
constant on $\sJ_A$-orbits, and let $P_A$ be the corresponding
orthogonal projection. The $P_A$ converge weakly to $P$.

\smallskip

\begin{prop}\label{orthproj}
Let  $A \supset B$ 
be finite sets of
non-archimedean places. Let $L$ be an $A$-invertible $K$-lattice, 
and $f\in H(B)$, the
restriction of $P_A f$ to the $\sJ_A$-orbit of $L$ is constant and
given by the formula
\begin{equation}\label{projformula}
P_A f|_{\sJ_A L} = \zeta_{K,A}(\beta)^{-1} \, \sum_{J\in \sJ_A}
\n(J)^{-\beta} \, f(J^{-1}L),
\end{equation}
where $\zeta_{K,A}(\beta)=\sum_{J\in \sJ_A} \n(J)^{-\beta}$.
\end{prop}

\proof By construction, the right hand side of the formula
\eqref{projformula} defines an element $f_A$ in $H_A \cap H(A)$. One checks,
using the quasi-invariance condition \eqref{qinvcond} on the
measure $\mu$, that $\langle f_A, g \rangle = \langle f, g\rangle$
for all $g\in H_A$, as in \cite{nes}.
\endproof

\smallskip

Let $L$ be an invertible $K$-lattice and $\chi$ a Grossencharakter
vanishing on the connected component of identity $D_K$. We define
$\chi(L)$ as $\chi(\rho s)$, for any representative $(\rho,s)$ of
$L$. This continues to make sense
when $L$ is an $A$-invertible $K$-lattice
and $A \supset B_\chi$ taking $\chi(\rho' s)$
where $\rho'\in \hat\O^*$ and $\rho'_v=\rho_v$ for all $v\in A$.

\smallskip

Finally we recall that to a Grossencharakter $\chi$
vanishing on the connected component of identity $D_K$ one
associates a Dirichlet character $\tilde{ \chi}$
defined for ideals 
  $J$  in $\sJ_{B_\chi^c}$, where $B_\chi^c$ is
the complement of $B_\chi$. More precisely, given 
 $J\in\sJ_{B_\chi^c}$, let $s_J$ be an id\`ele such that
$J=s_J \,\hat\O \cap K$ and $(s_J)_v=1$
for all places $v\in B_\chi$. One then define 
$\tilde{ \chi}(J)$ to be the value
$\chi(s_J)$. This is independent of the choice of such $s_J$.

\smallskip

\begin{prop}\label{PAfBchi}
Let $A\supset B\supset B_\chi$ and $L$ an $A$-invertible $K$-lattice.
The projection $P_A$ of \eqref{projformula} applied to the
function $f_{B,\chi}$ gives
\begin{equation}\label{projformulaB}
P_A f_{B,\chi}|_{\sJ_A L} = \frac{\chi(L)}{\zeta_{K,A}(\beta)} \,
\sum_{J\in \sJ_{A \backslash B}} \n(J)^{-\beta} \, \tilde{ \chi}(J)^{-1}.
\end{equation}
\end{prop}

\proof Among ideals in $\sJ_A$, those that have nontrivial
components on $B$ do not contribute to the sum \eqref{projformula}
computing $P_A f_{B,\chi}|_{\sJ_A L}$. It remains to show that
$f_{B,\chi}(J^{-1}L)=\chi(L)\tilde{ \chi}(J)^{-1}$, for $J\in \sJ_{A \backslash
B}$. Let $J=s_J \,\hat\O \cap K$ and $(s_J)_v=1$
for all places $v\in B$. Let $L$ be given by $(\rho,s)$, we have
$J^{-1}L$ given by $(\rho\,s_J , s\,s_J^{-1} )$ using 
Proposition \ref{actJL}. Thus,
for any choice of $\rho'\in \hat\O^*$ 
with $\rho'_v=(\rho\,s_J)_v$ for all $v\in B$ one has
$f_{B,\chi}(J^{-1}L)=\chi(\rho' s\,s_J^{-1})=\chi(\rho's)
\tilde{ \chi}(J)^{-1}$.
Note that $(s_J)_v=1$
for all places $v\in B$ thus the choice of $\rho'$
is governed by $\rho'_v=\rho_v$ for all $v\in B$.
Since $L$ is $A$-invertible and 
$A\supset B\supset B_\chi$ we get $\chi(\rho's)=\chi(L)$
for a suitable choice of $\rho'$.
\endproof

\smallskip

It then follows as in \cite{nes} that $P_A f_{B,\chi}$ tend weakly
to zero for $\chi$ nontrivial. The same argument gives an explicit
formula for the measure, obtained as a limit of the $P_A f_{B,1}$,
for the trivial character. In particular, the restriction of the
measure to $G_B$ is proportional to the Haar measure. Positivity
is ensured by the fact that we are taking a projective limit of
positive measures. This completes the proof of existence and
uniqueness of the KMS$_\beta$ state for $\beta\in (0,1]$.

\medskip
\subsection*{Open Questions}\hfill\medskip

Theorem \ref{KMSCM}
shows the existence of a $C^*$-dynamical system
$(\cA_K,\sigma_t)$ with
all the required properties for the interpretation
of the class field theory isomorphism in the CM case
in the framework of fabulous states. There is however
still one key feature of the BC-system that needs to be
obtained in this framework. It is the presentation of
the arithmetic subalgebra $\cA_{K,\Q}$ in terms of
generators and relations. This should be obtained along the
lines of \cite{CM} Section 6, Lemma 15 and Proposition 15,
and Section 9 Proposition 41. These suggest that the
relations will have coefficients in the Hilbert
modular field.

\smallskip

We only handled in this paper the CM-case \ie imaginary quadratic
fields, but many of the notions we introduced such as
that of a $K$-lattice should be extended to arbitrary
number fields $K$. Note in that respect that Proposition \ref{space1dKcomm}
indicates clearly that in general the space of
commensurability classes of $K$-lattices should be
identical to the space $\A_K/K^*$ of Ad\`ele classes introduced in
\cite{Co-zeta} for the spectral realization of zeros of
$L$-functions, with the slight nuance of non-zero
archimedan component. The scaling group which is used to
pass from the above ``dual system" to the analogue of
the BC system is given in the case $K=\Q$ by the group
$\R^*_+$ and in the case of imaginary quadratic fields
by the multiplicative group $\C^*$. It is thus natural
to expect in general that it will be given by the
connected component of identity $D_K$ in the group
$C_K$ of id\`ele
classes.

\end{document}